\documentclass[12pt,a4paper]{article}
\usepackage[a4paper,margin=1in]{geometry}
\usepackage[T1]{fontenc}
\usepackage{anyfontsize}
\usepackage{microtype}
\usepackage{amsmath,amssymb,amsthm,bm}
\usepackage{newtxtext,newtxmath}
\allowdisplaybreaks[4]
\usepackage{tikz}
\usetikzlibrary{positioning}

\usepackage[linesnumbered,ruled,vlined]{algorithm2e}
\usepackage{cite}
\usepackage[hidelinks]{hyperref}
\hypersetup{
  pdftitle={Level Decompositions for Symmetric Deformations of the Braid Arrangement},
  pdfauthor={Yanru Chen, Suijie Wang, Jinxing Yang, and Chengdong Zhao},
  pdfsubject={Region levels in symmetric deformations of the braid arrangement},
  pdfkeywords={hyperplane arrangements, braid arrangement, region level, labeled Dyck paths, Stirling convolution, Raney numbers}
}

\newtheorem{theorem}{Theorem}[section]
\newtheorem{lemma}[theorem]{Lemma}
\newtheorem{prop}[theorem]{Proposition}
\theoremstyle{definition}
\newtheorem{exam}[theorem]{Example}
\theoremstyle{plain}
\newtheorem*{claim}{Claim}

\newcommand{\A}{\mathcal{A}}
\numberwithin{equation}{section}
\begin{document}

\begin{center}
{\Large\bfseries
Level Decompositions for Symmetric Deformations of\\[2mm] the Braid Arrangement}\\[7pt]
\end{center}

\vskip 3mm

\begin{center}
 Yanru Chen$^{1}$, Suijie Wang$^{2}$, Jinxing Yang$^{3}$ and Chengdong Zhao$^{4}$\\[8pt]
 $^{1,2,3}$\,School of Mathematics\\
 Hunan University\\
 Changsha 410082, Hunan, P. R. China\\[12pt]

 $^{4}$
 School of Mathematics and Statistics\\
 Central South University\\
 Changsha 410083, Hunan, P. R. China\\[15pt]

 Emails: $^{1}$yanruchen@hnu.edu.cn, $^{2}$wangsuijie@hnu.edu.cn, \\~~~~~~~$^{3}$jxyangmath@hnu.edu.cn, $^{4}$cdzhao@csu.edu.cn\\[15pt]

\end{center}
\vskip 3mm

\begin{abstract}
	Let $A\subseteq\mathbb R_{\ge0}$ be finite and nonempty, and let
	$\mathfrak{A}^A=(\mathcal{A}_1^A,\mathcal{A}_2^A,\ldots)$ be the associated
	sequence of symmetric deformations of the braid arrangement. Denote by
	$r_\ell(\mathcal{A}_n^A)$ the number of its level-$\ell$ regions and by
	$F_\ell(\mathfrak{A}^A,x)$ the corresponding exponential generating
	function. We prove
    \[
    F_\ell(\mathfrak{A}^A,x)=\bigl(F_1(\mathfrak{A}^A,x)\bigr)^\ell .
    \]
	As a consequence, the characteristic polynomial has the binomial-basis
	expansion
    \[
    \chi(\mathcal{A}_n^A,t)=\sum_{\ell=1}^{n}(-1)^{n-\ell}\,r_\ell(\mathcal{A}_n^A)\binom{t}{\ell}.
    \]
	When $0\in A$ and $A^*=A\setminus\{0\}$ is nonempty, we refine a
	classical identity of Stanley level by level:
    \[
	F_\ell(\mathfrak{A}^{A^*},x)=F_\ell(\mathfrak{A}^{A},1-e^{-x}).
    \]
	Equivalently, the Catalan-type and semiorder-type level counts satisfy an
	unsigned Stirling convolution of the first kind. For the $m$-Catalan
	arrangement $\mathcal{A}_n^{[0,m]}$, we obtain
    \[
    r_\ell(\mathcal{A}_n^{[0,m]})=n!\,\operatorname{Ran}_{m+1,m\ell}(n-\ell),
    \]
	where $\operatorname{Ran}_{p,r}(q)$ is a Raney number. This realizes Raney
	numbers as refined region counts and answers a question of Deshpande,
	Menon, and Sarkar. The proofs use labeled Dyck paths, interval orders, and
	exponential sequences of arrangements. We also realize the inverse
	Fu--Wang--Zhu bijection for $m$-Catalan regions by tableaux.
\end{abstract}

\vskip 6pt

\noindent
\textbf{2020 Mathematics Subject Classification.} Primary 52C35; Secondary 05A15, 05A19, 06A07.
\\ [7pt]
\textbf{Keywords.} symmetric deformation of the braid arrangement; region level;
characteristic polynomial; Catalan-type arrangement; semiorder-type arrangement;
labeled Dyck path; Stirling convolution; Raney number.
\section{Introduction}

For a positive integer $q$, write $[q]=\{1,\ldots,q\}$ and
$[0,q]=\{0,1,\ldots,q\}$.
A \emph{deformation} of the braid arrangement in $\mathbb{R}^n$ is a
hyperplane arrangement defined by equations $x_i-x_j=a$, where $i\ne j$
and the permitted offsets form finite sets. Such deformations are central
objects in the enumerative theory of arrangements; see
\cite{Athanasiadis1996,P-S2000,Bernardi2018}. We call a deformation
\emph{symmetric} if
it is invariant under every permutation of the coordinates. It is then
determined by a nonempty finite set $A=\{a_1,\ldots,a_m\}$ of nonnegative real
numbers, which we order as $a_1>\cdots>a_m\ge0$. We write $\A_n^A$ for the
resulting arrangement, whose hyperplanes are
\begin{equation*}
	x_i-x_j=\pm a_1,\ldots,\pm a_m,\quad 1\le i<j\le n.
\end{equation*}
For a hyperplane arrangement $\A$ in $\mathbb{R}^n$, the connected
components of $\mathbb{R}^n\setminus\bigcup_{H\in\A}H$ are its
\emph{regions}. Write $R(\A)$ for the set of regions and
$r(\A)=|R(\A)|$. We refine the usual region enumeration by level. Following
Ehrenborg~\cite{Ehrenborg2019}, the level $\ell(X)$ of a subset
$X\subseteq\mathbb{R}^n$ is the minimum dimension of a linear subspace $W$
such that
\begin{equation}\label{def-level}
	X\subseteq B(W,r) = \bigl\{ {\bm x}\in\mathbb{R}^n : d({\bm x},W)\le r \bigr\}
\end{equation}
for some $r>0$, where $d$ is Euclidean distance. For regions of
arrangements, this agrees with the
\emph{degree of freedom} of Armstrong and Rhoades~\cite{Armstrong2012},
while $\ell(X)-1$ is the \emph{ideal dimension} of
Zaslavsky~\cite{Zaslavsky2003}. Thus level enumeration simultaneously
refines the total number of regions and the number of relatively bounded
regions. Let $R_\ell(\A)$ be the set of level-$\ell$ regions and put
$r_\ell(\A)=|R_\ell(\A)|$.

Every region of $\A_n^A$ is invariant under translation by the diagonal
line $\mathbb R(1,\ldots,1)$. Hence its level is at least $1$ and at most
$n$; in particular, $r_\ell(\A_n^A)=0$ unless $1\le\ell\le n$.

Fix a nonempty finite set $A$ of nonnegative real numbers. For each $n\ge 1$, let $\A_n^A$ be the corresponding symmetric deformation of the braid arrangement in $\mathbb{R}^n$, and write
\[
\mathfrak{A}^{A}=(\A_{1}^{A},\A_{2}^{A},\ldots)
\]
for the resulting sequence. The exponential generating function of the level-$\ell$ region counts is
\begin{equation}\label{def-function}
	F_{\ell}(\mathfrak{A}^A, x) =\sum_{n\ge 1}r_{\ell}(\A_{n}^A)\frac{x^n}{n!}.
\end{equation}

\begin{theorem}\label{main-1}
	For every finite nonempty $A\subseteq\mathbb R_{\ge0}$ and $\ell\ge 1$,
	we have
	\begin{equation*}
		F_{\ell}(\mathfrak{A}^A, x)=\big(F_{1}(\mathfrak{A}^A, x)\big)^\ell,
	\end{equation*}
	or equivalently, for every $n\ge1$,
	\begin{equation}\label{level-identity}
		r_{\ell}(\A_n^A)=\sum_{\substack{n_1+\cdots+n_{\ell}=n\\ n_{1},\ldots,n_\ell\ge 1}} \binom{n}{n_1,\ldots,n_\ell}r_1(\A_{n_1}^A)\cdots r_1(\A_{n_\ell}^A).
	\end{equation}
\end{theorem}
For a hyperplane arrangement $\A$, let $L(\A)$ be its \emph{intersection
poset}: it consists of the ambient space and all nonempty intersections of
nonempty subfamilies of $\A$, ordered by reverse inclusion. The
\emph{characteristic polynomial} of $\A$ is
\[
\chi(\A,t)=\sum_{X\in L(\A)}\mu(\mathbb{R}^n,X)\,t^{\dim X},
\]
where the M\"obius function $\mu$ is recursively defined by
$\mu(X,X)=1$ and
$\mu(X,Y)=-\sum_{X\le Z<Y}\mu(X,Z)$. The next result identifies the
binomial-basis coefficients with signed level counts.
\begin{theorem}\label{main-2}
	For every finite nonempty $A\subseteq\mathbb R_{\ge0}$ and $n\ge1$, we
	have
	\begin{equation*}
		\chi(\A_n^A,t) = \sum_{\ell=1}^{n}(-1)^{n-\ell}r_{\ell}(\A_n^A)\binom{t}{\ell}.
	\end{equation*}	
\end{theorem}

When $0\in A$ and $A^*=A\setminus\{0\}$ is nonempty, we call $\A_n^A$ the
\emph{Catalan-type} arrangement and $\A_n^{A^*}$ the
\emph{semiorder-type} arrangement. Stanley~\cite{Stanley1996} proved the
identity
\[
F(\mathfrak{A}^{A^*},x)=F(\mathfrak{A}^{A},1-e^{-x}),
\]
where $F(\mathfrak{A}^A,x)=\sum_{\ell\ge 1}F_{\ell}(\mathfrak{A}^A,x)$, which was later reproved by Postnikov and Stanley~\cite{P-S2000} via Burnside's lemma. We refine this identity level by level, using a tuple of labeled Dyck paths associated with each region.
\begin{theorem}\label{main-3} 
	If $0\in A$ and $A^*=A\setminus \{0\}$ is nonempty, then, for every
	$\ell\ge1$,
	\[
	F_{\ell}(\mathfrak{A}^{A^*},x)=F_{\ell}(\mathfrak{A}^{A},1-e^{-x}),
	\]
	or equivalently, for every $n\ge1$,
	\begin{equation}
		r_{\ell}(\A_n^A)=\sum_{k=1}^{n}c(n,k)r_{\ell}(\A_k^{A^*}),
		\label{Stirling-convolution}
	\end{equation}
	where $c(n,k)$ denotes the unsigned Stirling number of the first kind.
\end{theorem}

For the special case $A=[0,m]$, the arrangement $\A_n^A$
specializes to the classical \emph{$m$-Catalan arrangement}
$\A_n^{[0,m]}$. For integers $p,r\ge1$ and $q\ge0$, write
\[
\operatorname{Ran}_{p,r}(q)
=\frac{r}{pq+r}\binom{pq+r}{q}
\]
for the Raney number~\cite{Raney1960}. Thus
$\operatorname{Ran}_{m+1,1}(n)$ is the Fuss--Catalan number associated
with the $m$-Catalan arrangement. Stanley~\cite{Stanley1996} showed that
the number of regions of $\A_n^{[0,m]}$ is
$n!\,\operatorname{Ran}_{m+1,1}(n)$.
Deshpande, Menon, and Sarkar~\cite{DKW-2023} asked whether Raney numbers can
be realized through region counts of a natural family of arrangements. The
next theorem gives such a realization: after division by $n!$, the level
numbers are Raney numbers.
\begin{theorem}\label{main-4}
	For integers $m,n\ge1$ and $1\le \ell\le n$, the number of
	level-$\ell$ regions of
	$\A_n^{[0,m]}$ satisfies
	\[r_{\ell}(\A_n^{[0,m]})=n!\,\operatorname{Ran}_{m+1,m\ell}(n-\ell).\]
\end{theorem}

\medskip
\noindent
The four results form a single structural picture. Theorem~\ref{main-1}
decomposes a level-$\ell$ region into an ordered collection of $\ell$ level-one
pieces. Theorem~\ref{main-2} transfers this decomposition to the
characteristic polynomial through Stanley's theory of exponential sequences
of arrangements. Theorem~\ref{main-3} shows that Stanley's relation between
Catalan-type and semiorder-type arrangements respects every level, not only
the total region count. Finally, Theorem~\ref{main-4} specializes the
decomposition to the $m$-Catalan arrangement and identifies the normalized
level numbers with Raney numbers. The main combinatorial tools are labeled
Dyck paths and interval orders. The tableau reconstruction in
Section~\ref{sec-7} gives an inverse to the map of Fu, Wang, and
Zhu~\cite{F-W2021}, complementing the approach of Duarte and Guedes de
Oliveira~\cite{D-G2021}.

The principal results of this paper were first made public in the initial
version of \href{https://arxiv.org/abs/2410.10198}{arXiv:2410.10198},
posted on October 14, 2024. Since that version appeared, several works have
developed related aspects of region level. Chen, Fu, Wang, and Yang extended
the study to broader classes of deformed braid arrangements
\cite{ChenFuWangYang2026}. Southerland, Southern, and Zhou obtained a general
approach through centralization, including extensions to geometric
semilattices and alternative derivations of level-enumeration identities
\cite{SoutherlandSouthernZhou2025}. Zhang studied characteristic polynomials
of deformations of Coxeter arrangements through levels of regions
\cite{Zhang2026}. More recently, Chen, Fu, Liang, and Wang extended the
enumeration from regions to faces of exponential sequences of arrangements
and from characteristic polynomials to Whitney polynomials
\cite{ChenFuLiangWang2026}. This manuscript presents those 2024 results in
revised form, with expanded proofs and an updated account of subsequent
developments.

The paper is organized as follows. Section~\ref{sec-2} proves
Theorem~\ref{main-1} by decomposing labeled Dyck paths, treating
Catalan-type and semiorder-type arrangements separately.
Section~\ref{sec-3} derives Theorem~\ref{main-2} from Stanley's theory of
exponential sequences of arrangements. Section~\ref{sec-4} proves
Theorem~\ref{main-3} by a bijection involving cycle decompositions of
permutations, and Section~\ref{sec-6} proves Theorem~\ref{main-4} by
generating functions and binomial identities. Appendix~\ref{sec-7}
presents the tableau realization of the inverse Fu--Wang--Zhu bijection.

\section{Level decompositions via labeled Dyck paths}\label{sec-2}
We prove Theorem~\ref{main-1}. Expanding the $\ell$-th power yields
\[
\begin{aligned}
\big(F_1(\mathfrak{A}^A,x)\big)^\ell
&=\Bigl(\sum_{n\ge 1}r_{1}(\A_n^A)\frac{x^n}{n!}\Bigr)^{\!\ell}\\
&=\sum_{n\ge 1}\frac{x^n}{n!}
\sum_{\substack{n_1,\ldots,n_\ell\ge 1\\ n_{1}+\cdots+n_{\ell}=n}}
\binom{n}{n_1,\ldots,n_\ell}\prod_{j=1}^\ell r_1(\A_{n_j}^A).
\end{aligned}
\]
Let $R_{\ell}(\A_n^A)$ be the set of level-$\ell$ regions of $\A_n^A$ and $r_\ell(\A_n^A)=|R_\ell(\A_n^A)|$. It therefore suffices to construct a bijection
\[
R_{\ell}(\A_n^A)\;\longrightarrow\;\bigsqcup_{\substack{ B_{1},\ldots,B_\ell\neq \emptyset \\ B_1\sqcup\cdots\sqcup B_{\ell}=[n]}}
\{(B_1,\ldots,B_\ell)\}\times R_1(\A_{|B_1|}^A)\times\cdots\times R_1(\A_{|B_\ell|}^A),
\]
where $\sqcup$ denotes disjoint union and $[n]=\{1,\ldots,n\}$, which is equivalent to the identity~\eqref{level-identity}. Our main tool for this purpose is the \emph{labeled Dyck path}, introduced below.

A \emph{Dyck path} $D$ of length $2n$ is a lattice path from $(0,n)$ to
$(n,0)$ with $n$ East and $n$ South steps that stays weakly above the
diagonal $x+y=n$. It is \emph{prime} if it meets the diagonal only at its
endpoints. Every Dyck path decomposes uniquely into concatenated prime
paths, called its \emph{prime components}; let $c(D)$ be their number. A
\emph{labeled Dyck path} $D_\pi$ is a Dyck path together with
$\pi\in\mathfrak{S}_n$: the $i$th East step from the left and the $i$th
South step from the top both carry the label $\pi(i)$.

\subsection{Catalan-type arrangements}
If $A=\{0\}$, then $\A_n^A$ is the braid arrangement: all its $n!$
regions have level $n$. Hence $F_1(\mathfrak{A}^A,x)=x$ and
$F_\ell(\mathfrak{A}^A,x)=x^\ell$, so Theorem~\ref{main-1} holds directly.
For the rest of this subsection, assume $0\in A$ and write
$A=\{0,a_1,\ldots,a_m\}$ with $a_1>\cdots>a_m>0$. The Catalan-type arrangement $\A_n^A$ consists of the hyperplanes
\[ x_i-x_j=0,\;\pm a_1,\ldots,\pm a_m,\qquad 1\le i<j\le n.\]
Observe that $\A_n^A$ contains the braid arrangement $\mathcal{B}_n=\{x_i-x_j=0:1\le i<j\le n\}$, whose regions are the open cones
\[
C_\pi:\; x_{\pi(1)}>x_{\pi(2)}>\cdots>x_{\pi(n)}
\]
for $\pi\in\mathfrak{S}_n$. Each cone $C_\pi$ is subdivided in exactly the same way by the additional hyperplanes of $\A_n^A$. We refer to $C_e$ as the {\it fundamental chamber}.
The symmetric group $\mathfrak{S}_n$ acts on $\mathbb{R}^n$ by the standard
left coordinate action
\[
(\pi\cdot{\bm x})_i=x_{\pi^{-1}(i)}.
\]
This induces an action on the regions of $\A_n^A$ via
  \[
 \Delta\longmapsto \pi\cdot\Delta
 =\bigl\{\pi\cdot{\bm x}:{\bm x}\in\Delta\bigr\},
 \]
which maps regions inside the fundamental chamber $C_e$ bijectively onto
regions inside $C_\pi$, where $e$ is the identity permutation. Since
$\pi$ is an isometry, $\Delta\subseteq B(W,r)$ if and only if
$\pi\cdot\Delta\subseteq B(\pi\cdot W,r)$; hence
$\ell(\Delta)=\ell(\pi\cdot\Delta)$ and the bijection preserves levels.

Thus exactly $\frac{1}{n!}$ of the level-$\ell$ regions of $\A_n^A$ lie in $C_e$. Let $R_{e,\ell}(\A_n^A)$ be the set of such regions and $r_{e,\ell}(\A_n^A)=|R_{e,\ell}(\A_n^A)|$. Then $r_\ell(\A_n^A)=n!\,r_{e,\ell}(\A_n^A)$, and identity~\eqref{level-identity} reduces to
\[
r_{e,\ell}(\A_n^A)=\sum_{\substack{n_1,\ldots,n_\ell\ge 1\\ n_1+\cdots+n_\ell=n}} r_{e,1}(\A_{n_1}^A)\cdots r_{e,1}(\A_{n_\ell}^A).
\]
To prove this identity, we associate to each region $\Delta\subseteq C_e$ a sign matrix $M_1(\Delta)$ and a Dyck path $D_1(\Delta)$. For any ${\bm x}=(x_1,\ldots,x_n)\in\Delta$, define
\begin{align}\label{catalan-M1}
M_1(\Delta) =\bigl(\operatorname{sgn}(x_i-x_j-a_1)\bigr)_{i,j\in[n]},
\end{align}
where $\operatorname{sgn}$ assigns to a nonzero real number its sign $+$ or $-$. (The subscript $1$ indicates that the matrix is built from the largest parameter $a_1$.) The matrix is well defined because $x_i-x_j-a_1\neq0$ for every $i,j$, and its entries remain constant as ${\bm x}$ varies within $\Delta$. Since $\Delta\subseteq C_e$, we have $x_1>\cdots>x_n$, which forces all $+$ entries to lie above the diagonal and cluster in the northeast corner of $M_1(\Delta)$. The southwest boundary separating the $+$ and $-$ regions therefore forms a Dyck path of length $2n$, which we denote by $D_1(\Delta)$: entries are $+$ above $D_1(\Delta)$ and $-$ below it. We illustrate the construction with an example.
\begin{exam} Let $A=\{0,1,2,3\}$ and let $\Delta$ be a region of $\A_5^A$ in the fundamental chamber $C_e$ defined by
\begin{align*}
\Delta = \left\{ (x_{1}, \ldots, x_{5}) \in \mathbb{R}^{5} \; \middle| \; 
\begin{aligned}
& x_{1} > x_{2} > x_{3} > x_{4} > x_{5}, \\
& 2> x_{1} - x_{2} > 1, \\
& x_{2} - x_{3}>3, \\
& 3>x_3 -x_4>2,\\
& 3>x_{4} - x_{5}>2.
\end{aligned}
\right\}
\end{align*}
The sign matrix $M_1(\Delta)$ determines the Dyck path $D_1(\Delta)$,
shown in red in Figure~\ref{fig:cat-1}.
\begin{figure}[htbp]
\centering 

\tikzset{global scale/.style={
scale=#1,
every node/.append style={scale=#1}}} 
\begin{tikzpicture}[global scale = 0.9]
\begin{scope}[scale=0.8, shift={(0,15)}] 
\draw (0,0) grid (5,5);

\draw[red,line width=0.5mm ] (0,5) -- (2,5) -- (2,3) -- (3,3) -- (4,3) -- (4,2) -- (5,2) -- (5,0);

\node at (-0.7,5) {$(0,5)$};
\node at (5,-0.4) {$(5,0)$};
\node at (2.5,-1) {$D_1(\Delta)$};

\node at (-0.2,-0.2) {$O$};
\draw[->, thick] (0,0) -- (5.8,0) node[below right] {$x$};
\draw[->, thick] (0,0) -- (0,5.8) node[above left] {$y$};

\node[] at (-6,2.5)
{
	$
	\renewcommand{\arraystretch}{1.6}
	\setlength{\arraycolsep}{9pt}
	\begin{pmatrix}
		- & - & + & + &+\\
		- & - & + & + &+ \\
		- & - & - & - &+\\
		- & - & - & - &-\\
		- & - & - & - &- \\
	\end{pmatrix}
	$
};
\node at (-6,-1) {$M_1(\Delta)=(\operatorname{sgn}(x_{i} - x_{j} - 3))_{i,j}$};
\end{scope}
\end{tikzpicture}
\caption{The sign matrix $M_1(\Delta)$ and Dyck path $D_1(\Delta)$}
\label{fig:cat-1}
\end{figure}

By Proposition~\ref{components-1}, the number of prime components of $D_1(\Delta)$ is $c(D_1(\Delta))=2$, which implies that the level of $\Delta$ is $\ell(\Delta)=2$.
\end{exam}
\begin{prop}\label{components-1}
Let $\Delta$ be a region of $\A_n^A$ contained in the fundamental chamber
$C_e$. Then the level of $\Delta$ is the number of prime components of
$D_1(\Delta)$; that is, $\ell(\Delta)=c(D_1(\Delta))$.
\end{prop}
\begin{proof}
Let $\Delta$ be a region of $\A_n^A$ inside $C_e$, and suppose
$D_1(\Delta)$ has $\ell$ prime components intersecting the diagonal at
$(i_0,n-i_0),\ldots,(i_\ell,n-i_\ell)$, where
$0=i_0<\cdots<i_\ell=n$. For $1\le i<n$, the path returns to the diagonal
at $(i,n-i)$ exactly when $(M_1(\Delta))_{i,i+1}=+$, equivalently when
$x_i-x_{i+1}>a_1$. In that case
$x_p-x_q>a_1$ for all $p\le i<q$. Thus, for every
$(x_1,\ldots,x_n)\in\Delta$,
\begin{equation}\label{linear-order}
\begin{cases}
    x_i-x_{i+1}>a_1, &\text{if }i\in\{i_1,\ldots,i_{\ell-1}\},\\[2pt]
    x_i-x_{i+1}<a_1, &\text{otherwise}.
\end{cases}
\end{equation}
Define the $\ell$-dimensional subspace
\begin{align*}
    W = \bigl\{ (y_1,\ldots,y_n)\in\mathbb{R}^n : y_i=y_j\;\text{whenever}\; i_{k-1}< i,j\le i_k \text{ for some } k \bigr\},
\end{align*}
so that $\dim W=\ell$. For any ${\bm x}=(x_1,\ldots,x_n)\in\Delta$, choose ${\bm z}=(z_1,\ldots,z_n)\in W$ by setting
\begin{equation*}
    \begin{cases}
        z_1 = \cdots = z_{i_1} = x_{i_1},\\[3pt]
        z_{i_1+1} = \cdots = z_{i_2} = x_{i_2},\\[3pt]
        \qquad\vdots\\[3pt]
        z_{i_{\ell-1}+1} = \cdots = z_n = x_n.
    \end{cases}
\end{equation*}
Within each block, every adjacent difference is less than $a_1$ by
\eqref{linear-order}. Hence
$0\le x_t-x_{i_k}\le(i_k-t)a_1\le(n-1)a_1$, and
\[
d({\bm x},{\bm z})
=\left(\sum_{k=1}^\ell\sum_{t=i_{k-1}+1}^{i_k}(x_t-x_{i_k})^2\right)^{1/2}
\le\sqrt n\,(n-1)a_1.
\]
Thus $\Delta\subseteq B(W,1+\sqrt n\,(n-1)a_1)$, which gives
$\ell(\Delta)\le \dim W=\ell$.

For the reverse inequality, let $C\subseteq W$ be the cone
\[
C=\{(y_1,\ldots,y_n)\in W: y_1\ge\cdots\ge y_n\}.
\]
The cone $C$ spans $W$. Fix ${\bm x}_0\in\Delta$ and choose a linear
subspace $U$ of dimension $\ell(\Delta)$ and $r>0$ such that
$\Delta\subseteq B(U,r)$. Since ${\bm x}_0+C\subseteq\Delta$, for every
${\bm c}\in C$ we have
\[
d({\bm c},U)\le \|{\bm x}_0\|+d({\bm x}_0+{\bm c},U)
\le \|{\bm x}_0\|+r.
\]
If ${\bm c}\in C$, the same bound applies to $t{\bm c}$ for every
$t>0$; dividing by $t$ and letting $t\to\infty$ gives
${\bm c}\in U$. Thus $C\subseteq U$, and hence
$W=\operatorname{span}C\subseteq U$. Therefore
$\ell\le\dim U=\ell(\Delta)$.
\end{proof}

Write $e$ for the identity permutation in each $\mathfrak{S}_n$. We define a map
\begin{equation}\label{catalan-map}
    \varphi_{\ell}\colon R_{e, \ell}(\A_n^A) \longrightarrow 
    \bigsqcup_{\substack{n_1,\ldots,n_\ell\ge 1\\ n_1+\cdots+n_\ell=n}}R_{e,1}(\A_{n_1}^A)\times\cdots\times R_{e,1}(\A_{n_\ell}^A)
\end{equation}
as follows. For $\Delta\in R_{e,\ell}(\A_n^A)$, Proposition~\ref{components-1} guarantees that $D_1(\Delta)$ has $\ell$ prime components. Let $(i_0,n-i_0),(i_1,n-i_1),\ldots,(i_\ell,n-i_\ell)$ be the intersection points of $D_1(\Delta)$ with $x+y=n$, where $0=i_0<\cdots<i_\ell=n$, and set
\[
n_1=i_1,\quad n_2=i_2-i_1,\quad\ldots,\quad n_\ell=n-i_{\ell-1}.
\]
Take ${\bm x}=(x_1,\ldots,x_n)\in\Delta$. For each $1\le k\le \ell$, the subvector $(x_{i_{k-1}+1},\ldots,x_{i_k})$ belongs to a uniquely determined region $\Delta_k$ of $\A_{n_k}^A$. We set
\[
\varphi_\ell(\Delta)=(\Delta_1,\Delta_2,\ldots,\Delta_\ell).
\]

\begin{theorem}
	For every integer $\ell$ with $1\le\ell\le n$, the map $\varphi_{\ell}$ defined
	in~\eqref{catalan-map} is bijective.
\end{theorem}
\begin{proof}
We first check that $\varphi_\ell$ is well defined.
\begin{enumerate}
    \item[(1)] \emph{Independence of the choice of ${\bm x}$.} Suppose
    ${\bm x},{\bm y}\in\Delta$ yield distinct images. Then, for some $k$,
    the corresponding subvectors lie in different regions of
    $\A_{n_k}^A$. Hence there are $1\le p<q\le n_k$ and $1\le s\le m$
    such that
    \[
    x_{i_{k-1}+p}-x_{i_{k-1}+q}>a_s,
    \qquad
    y_{i_{k-1}+p}-y_{i_{k-1}+q}<a_s,
    \]
    contradicting that ${\bm x}$ and ${\bm y}$ lie in the same region
    $\Delta$.

    \item[(2)] \emph{Each $\Delta_k$ has level $1$.} By the construction of $\Delta_k$, we have that the sign matrix $M_1(\Delta_k)$ is precisely the $n_k\times n_k$ principal submatrix of $M_1(\Delta)$ indexed by $\{i_{k-1}+1,\ldots,i_k\}$. Consequently, $D_1(\Delta_k)$ coincides with the $k$-th prime component of \allowbreak $D_1(\Delta)$, so $\ell(\Delta_k)=1$.
\end{enumerate}

We prove injectivity. Assume $\Delta\neq\Delta'$ in $R_{e,\ell}(\A_n^A)$ satisfy $\varphi_\ell(\Delta)=\varphi_\ell(\Delta')=(\Delta_1,\ldots,\allowbreak \Delta_\ell)$ with $\Delta_k\in R_{e,1}(\A_{n_k}^A)$. Then $D_1(\Delta)$ and $D_1(\Delta')$ share the same intersection points $(i_0,n-i_0),\ldots,(i_\ell,n-i_\ell)$ with $x+y=n$, where $i_k-i_{k-1}=n_k$. Since $\Delta\neq\Delta'$, there exist $i<j$ and $1\le s\le m$ such that for all ${\bm x}\in\Delta$, ${\bm x}'\in\Delta'$,
\begin{equation}\label{contradiction1}
x_i-x_j>a_s\quad\text{and}\quad x'_i-x'_j<a_s.
\end{equation}
If $i,j$ lie in the same block (i.e., $i_{k-1}<i<j\le i_k$), then the corresponding subvectors belong to the same region $\Delta_k$, contradicting~\eqref{contradiction1}. If $i\le i_k<j$, then~\eqref{linear-order} forces $x_i-x_j>a_1$ and $x'_i-x'_j>a_1$, also contradicting~\eqref{contradiction1}. Thus no such pair exists, a contradiction.

Finally, we establish surjectivity. Fix $n_1+\cdots+n_\ell=n$ with each $n_k\ge 1$, and let $(\Delta_1,\ldots,\Delta_\ell)$ be an $\ell$-tuple with $\Delta_k\in R_{e,1}(\A_{n_k}^A)$. Choose ${\bm x}^{(k)}=(x^{(k)}_1,\ldots,x^{(k)}_{n_k})\in\Delta_k$. Since $(x^{(k)}_1,\ldots,x^{(k)}_{n_k})+t(1,\ldots,1)\in \Delta_k$ holds for any $t\in\mathbb{R}$, we may arrange that
\[
x^{(p)}_{n_p}-x^{(p+1)}_1>a_1\qquad(p=1,\ldots,\ell-1).
\]
Form the concatenated vector
\[
{\bm x}=(x^{(1)}_1,\ldots,x^{(1)}_{n_1},x^{(2)}_1,\ldots,x^{(2)}_{n_2},\ldots,x^{(\ell)}_1,\ldots,x^{(\ell)}_{n_\ell}),
\]
and let $\Delta$ be the region of $\A_n^A$ containing ${\bm x}$. The
coordinate inequalities inside the blocks and the chosen translation gaps
give $\Delta\subseteq C_e$. Since each $D_1(\Delta_k)$ is prime, we have $\operatorname{sgn}(x^{(k)}_i-x^{(k)}_{i+1}-a_1)=-$, while the construction ensures $\operatorname{sgn}(x^{(k)}_{n_k}-x^{(k+1)}_1-a_1)=+$. Hence $(M_1(\Delta))_{i,i+1}=+$ exactly when $i=n_1+\cdots+n_j$ for $1\le j<\ell$, which means $D_1(\Delta)$ consists of $\ell$ prime components of lengths $n_1,\ldots,n_\ell$. Therefore $\Delta\in R_{e,\ell}(\A_n^A)$ and $\varphi_\ell(\Delta)=(\Delta_1,\ldots,\Delta_\ell)$.
\end{proof}

\subsection{Semiorder-type arrangements}
Throughout this subsection, assume $0\notin A$ and write $A=\{ a_1,\ldots, a_m\}$ with $a_1>\cdots>a_m>0$. The semiorder-type arrangement $\A_n^A$ is given by the hyperplanes
\[ x_i-x_j=\pm a_1,\ldots, \pm a_m, \quad  1\le i< j\le n.\]
Stanley \cite{Stanley1996} applied interval orders to the study of hyperplane arrangements, particularly to semiorder-type arrangements. We explain the interval order concept, our labeled Dyck path model, and the connection between them.

Recall that for a finite set $P = \{I_1, \ldots, I_n\}$ consisting of closed intervals \( I_i = [a_i, b_i] \) with \( a_i, b_i \in \mathbb{R} \) and \( a_i < b_i \), we define a partial order on \(P\) by \( I_i < I_j \) if \( b_i < a_j \). We allow duplicated intervals in $P$, treating them as distinct. A poset isomorphic to such \(P\) is an {\it interval order}. When all intervals $I_i$ have the same length, $P$ is called a semiorder.

We associate a semiorder on $[n]$ with each region $\Omega$ of $\A_n^A$.
For ${\bm x}=(x_1,\ldots,x_n)\in\Omega$, assign to $i$ the interval
$[x_i-a_1,x_i]$ and define
\[
i\prec_{P_1(\Omega)}j
\quad\Longleftrightarrow\quad
x_i<x_j-a_1.
\]
Because all intervals have length $a_1$, this is a semiorder. It is
independent of ${\bm x}$: every comparison is determined by the sign of
$x_j-x_i-a_1$, which is constant on $\Omega$. Let $G[P_1(\Omega)]$ be its
incomparability graph, with vertex set $[n]$ and an edge $\{i,j\}$ exactly
when $i$ and $j$ are incomparable.

\begin{exam}\label{ex:semiorder-running}
For the region $\Omega$ of $\A_5^{[3]}$ below,
Figure~\ref{fig:poset-io} shows the semiorder $P_1(\Omega)$ and its
incomparability graph $G[P_1(\Omega)]$.
\[
\Omega=\left\{(x_{1},x_2,x_3,x_4,x_{5})\in \mathbb{R}^{5}\; \middle| \;
\begin{aligned}
&-1<x_{1}-x_{2}<1,\\
&1<x_{3}-x_{1}<2,\\ 
&2<x_{3}-x_{2}<3,\\
&3<x_{4}-x_{3}, \,\, x_{5}-x_{3}>3,\\
&-1<x_{4}-x_{5}<1.
\end{aligned}
\right\}
\]
\begin{figure}
\centering
\begin{tikzpicture}[node distance=1.5cm, every node/.style={circle,draw,fill=black,inner sep=2pt}, line width=0.3mm]

\begin{scope}[shift={(0,0)}]
\node[label=left:4] (1r) at (1,2) {};
\node[label=below:3] (2r) at (0,0) {};
\node[label=below:2] (3r) at (4,0) {};
\node[label=below:1] (4r) at (2,0) {};
\node[label=right:5] (5r) at (3,2) {};
\draw (1r) -- (2r) -- (5r);
\draw (1r) -- (4r) -- (5r);
\draw (1r) -- (3r) -- (5r);
\node[shape=rectangle, draw=none, fill=none, inner sep=0pt, outer sep=0pt] at (2,-1.2) {$P_1(\Omega)$};
\end{scope}

\begin{scope}[shift={(6,0)}]
	\node[label=left:4] (1r) at (1,2) {};
	\node[label=below:3] (2r) at (0,0) {};
	\node[label=below:2] (3r) at (4,0) {};
	\node[label=below:1] (4r) at (2,0) {};
	\node[label=right:5] (5r) at (3,2) {};
	\draw (2r) -- (4r) -- (3r);
	\draw (1r) -- (5r);
	\draw (2r) to[bend left] (3r);
	\node[shape=rectangle, draw=none, fill=none, inner sep=0pt, outer sep=0pt] at (2,-1.2) {$G[P_1(\Omega)]$};
\end{scope}
\end{tikzpicture}
\caption{The corresponding semiorder $P_1(\Omega)$ and incomparability graph $G[P_1(\Omega)]$}
\label{fig:poset-io}
\end{figure}
\end{exam}

\begin{prop}
	For any region $\Omega$ of $\A_n^A$, let $G[P_1(\Omega)]$ be the incomparability graph of the semiorder $P_1(\Omega)$. Then
	\[
	\ell(\Omega)=\kappa(G[P_1(\Omega)]),
	\]
	where $\kappa(G[P_1(\Omega)])$ denotes the number of connected components
	of $G[P_1(\Omega)]$.
\end{prop}
\begin{proof}
Let $G_1,\ldots,G_\ell$ be the connected components of
$G[P_1(\Omega)]$. If $i$ and $j$ are joined by an edge, then
$|x_i-x_j|<a_1$ throughout $\Omega$. A path inside a component has at most
$n-1$ edges and therefore shows that
\[
|x_i-x_j|<na_1
\qquad(i,j\in G_p).
\]

We use the standard fact that the connected components of the
incomparability graph of a poset are totally ordered. Indeed, elements in
distinct components are comparable. Moreover, the direction of comparison
is constant on each pair of components: if $u$ and $u'$ are incomparable
and $u\prec v$, then $v\prec u'$ would imply $u\prec u'$, a contradiction;
the assertion follows by propagating along incomparability paths. Transitivity
then orders the components. Relabel them so that
\[
x_i-x_j>a_1
\qquad(i\in G_p,\ j\in G_q, p<q)
\]
throughout $\Omega$.

Let $W$ be the $\ell$-dimensional subspace of vectors that are constant on each
$G_p$. Replacing the coordinates in each component by one representative
shows that every point of $\Omega$ lies within distance
$\sqrt n\,(n-1)a_1$ of $W$. Hence $\ell(\Omega)\le \ell$.

For the reverse inequality, fix ${\bm x}_0\in\Omega$ and let $C\subseteq W$
be the full-dimensional cone whose component values weakly decrease in the
order $G_1,\ldots,G_\ell$. Adding a vector of $C$ preserves every defining
inequality of $\Omega$, so ${\bm x}_0+C\subseteq\Omega$. If
$\Omega\subseteq B(U,r)$ for a linear subspace $U$, the scaling argument
used in Proposition~\ref{components-1} gives $C\subseteq U$. Since
$\operatorname{span}C=W$, we obtain $\ell=\dim W\le\dim U$. Therefore
$\ell\le\ell(\Omega)$.
\end{proof}

As in the Catalan-type case, we encode semiorder-type regions by labeled
Dyck paths. Fix a region $\Omega$ of $\A_n^A$ and a point
${\bm x}=(x_1,\ldots,x_n)\in\Omega$. Its \emph{associated permutation} is
the unique $\pi\in\mathfrak{S}_n$ such that
$x_{\pi(1)}\ge\cdots\ge x_{\pi(n)}$, with smaller labels first among
equal coordinates. Define
\begin{align}\label{semiorder-M1}
M_1({\bm x}) = \bigl(\operatorname{sgn}(x_{\pi(i)}-x_{\pi(j)}-a_1)\bigr)_{i,j\in[n]},
\end{align}
where $\operatorname{sgn}$ takes the values $+$, $-$, or $0$. (Within a region, $0$ never occurs.) Because $x_{\pi(1)}\ge\cdots\ge x_{\pi(n)}$, all $+$ entries lie strictly above the diagonal and are concentrated in the northeast corner of $M_1({\bm x})$. The boundary separating $+$ and $-$ forms a labeled Dyck path $D_{\pi,1}(\bm x)$ of length $2n$. We illustrate the construction with the following example.
\begin{exam}\label{ex2}
	Continue with the region $\Omega$ from
	Example~\ref{ex:semiorder-running}.
	Take
	\[
	{\bm x}=(3.0,2.2,4.5,8.0,8.0)\in\Omega.
	\]
	Its associated permutation is $\pi=45312$; the labeled Dyck path
	$D_{\pi,1}({\bm x})$ is shown in Figure~\ref{fig-semiorder-1}.
	
	\begin{figure}[htbp]
		\centering  
		\tikzset{global scale/.style={
				scale=#1,
				every node/.append style={scale=#1}
			}
		}  
		\begin{tikzpicture}[global scale = 1]
			\begin{scope}[scale=0.8, shift={(-13,6.5)}]
				\draw (0,0) grid (5,5);
				\draw[red,line width=0.5mm] (0,5) -- (2,5) -- (2,3) -- (5,3) -- (5,0);
				\node at (-0.2,-0.2) {$O$};
				\draw[->, thick] (0,0) -- (5.8,0) node[below right] {$x$};
				\draw[->, thick] (0,0) -- (0,5.8) node[above left] {$y$};
				\node at (0.5,5.3) {4};
				\node at (1.5,5.3) {5}; 
				\node at (1.75,4.5) {4};
				\node at (2.5,2.65) {3};
				\node at (1.75,3.5) {5};
				\node at (5.3,2.5) {3};
				\node at (3.5,2.65) {1}; 
				\node at (5.3,1.5) {1};
				\node at (4.5,2.65) {2};
				\node at (5.3,0.5) {2}; 
				\node at (2.5,-0.6) {$D_{\pi,1}(\bm x)$};
				\node[] at (-6,2.7)
				{
					$
					\renewcommand{\arraystretch}{1.6}
					\setlength{\arraycolsep}{9pt}
					\begin{pmatrix}
						- & - & + & + &+\\
						- & - & + & + &+ \\
						- & - & - & - &-\\
						- & - & - & - &-\\
						- & - & - & - &- \\
					\end{pmatrix}
					$
				};
				\node at (-6,-0.6) {$M_1(\bm x)=(\operatorname{sgn}(x_{\pi(i)} - x_{\pi(j)} - 3))_{i,j}$};
			\end{scope}  
		\end{tikzpicture}
		\caption{The sign matrix $M_1(\bm x)$ and labeled Dyck path $D_{\pi,1}(\bm x)$}
		\label{fig-semiorder-1}
	\end{figure}
\end{exam}

\begin{lemma}\label{matrix-independent}
	Let $\Omega$ be a region of $\A_n^A$. Then the matrix $M_{1}({\bm x})$ is independent of the choice of ${\bm x}\in \Omega$.
\end{lemma}

\begin{proof}
Let ${\bm x},{\bm y}\in\Omega$. Since every region of a real
hyperplane arrangement is an open convex polyhedron, they can be joined by
a path in $\Omega$. After an arbitrarily small perturbation inside
$\Omega$, we may assume that along this path at most two coordinates are
equal at any time. Between such equality times, the associated permutation
is constant, and all signs defining $M_1$ are constant because the path
does not leave $\Omega$.

It remains to examine a time at which two coordinates, say those labeled
$u$ and $v$, cross. Immediately before and after the crossing, the
associated permutations differ by the adjacent transposition of $u$ and
$v$. At the crossing $x_u=x_v$. Hence, for every label $w$,
\[
\begin{aligned}
\operatorname{sgn}(x_u-x_w-a_1)
 &=\operatorname{sgn}(x_v-x_w-a_1),\\
\operatorname{sgn}(x_w-x_u-a_1)
 &=\operatorname{sgn}(x_w-x_v-a_1).
\end{aligned}
\]
The two rows and the two columns exchanged in $M_1$ are therefore
identical at the crossing; the entries involving only $u$ and $v$ are all
negative because $a_1>0$. Swapping $u$ and $v$ consequently leaves the
unlabeled sign matrix unchanged. Repeating this argument at every crossing
shows that $M_1({\bm x})=M_1({\bm y})$.
\end{proof}
Lemma~\ref{matrix-independent} guarantees that all points in the same region yield the same unlabeled Dyck path. Given a region $\Omega\in R(\A_n^A)$ and $\bm x\in \Omega$, suppose $D_{\pi,1}(\bm x)$ has $\ell$ prime components. Define the ordered partition $\lambda(D_{\pi,1}(\bm x))=(B_1,\ldots,B_\ell)$ of $[n]$ where $B_i$ is the set of indices corresponding to the $i$-th prime component of $D_{\pi,1}(\bm x)$.

\begin{lemma}\label{partition-independent}
	Let $\Omega$ be a region of $\A_n^A$. Then the ordered partition $\lambda(D_{\pi,1}(\bm x))$ is independent of the choice of $\bm x\in\Omega$.
\end{lemma}
\begin{proof}
	Given $\bm x\in\Omega$, let $\pi$ be its associated permutation and $\lambda(D_{\pi,1}(\bm x))=(B_1,\ldots,B_\ell)$ with $|B_k|=n_k$. By construction of $D_{\pi,1}(\bm x)$, $\operatorname{sgn}(x_{\pi(i)}-x_{\pi(i+1)}-a_1)=+$ in $M_1(\bm x)$ if and only if $i=n_1+\cdots+n_k$ for $1\le k<\ell$, i.e.,
	\begin{equation*}
	\begin{cases}
		x_{\pi(i)}-x_{\pi(i+1)}-a_1>0, \;\text{ $i=n_1, n_1+n_2,\ldots,n_1+\cdots+n_{\ell-1}$};\\
		x_{\pi(i)}-x_{\pi(i+1)}-a_1<0, \;\text{ otherwise}.
	\end{cases}
	\end{equation*}
	Suppose there exists $\bm y\in\Omega$ with associated permutation $\pi'$ such that $\lambda(D_{\pi',1}(\bm y))$ differs from $\lambda(D_{\pi,1}(\bm x))$. By Lemma~\ref{matrix-independent}, both labeled paths have $\ell$ prime components with the same component sizes. Write
	$\lambda(D_{\pi',1}(\bm y))=(B'_1,\ldots,B'_\ell)$ and let $p$ be the
	first index for which $B_p\ne B'_p$. Choose
	$i_0\in B_p\setminus B'_p$ and $j_0\in B'_p\setminus B_p$. Since all
	earlier blocks agree and corresponding blocks have equal size, $i_0$
	belongs to a later block of the primed partition and $j_0$ to a later
	block of the unprimed partition. Therefore
	\[
	x_{i_0}-x_{j_0}>a_1 \text{ and } y_{j_0}-y_{i_0}>a_1,
	\]
	contradicting that $\bm x,\bm y$ lie in the same region.
\end{proof}
\begin{prop}
	Let $\Omega$ be a region of $\A_n^A$ and $\bm x$ an arbitrary point in $\Omega$. Then $\ell(\Omega)=c(D_{\pi,1}(\bm x))$.
\end{prop}
\begin{proof}
	Lemma~\ref{partition-independent} gives an ordered partition
	$(B_1,\ldots,B_\ell)$ of $[n]$. Elements in a common block are connected by
	a chain of pairs whose coordinate differences have absolute value less
	than $a_1$; hence $|x_i-x_j|<na_1$ for $i,j\in B_p$. For
	$i\in B_p$ and $j\in B_q$ with $p<q$, one has $x_i-x_j>a_1$.
	The upper- and lower-bound arguments in the preceding proposition, with
	the $B_p$ in place of the components $G_p$, now give
	$\ell(\Omega)=\ell$.
\end{proof}
Now we define a map
\begin{equation}\label{semiorder-map}
\phi_\ell\colon R_\ell(\A_n^A)\rightarrow
\bigsqcup_{\substack{B_{1},\ldots,B_\ell\neq \emptyset \\ B_1\sqcup \cdots\sqcup B_{\ell}=[n]}}
\{(B_1,\ldots,B_\ell)\}\times R_1(\A_{|B_1|}^A)\times\cdots\times R_1(\A_{|B_\ell|}^A),
\end{equation}
as follows. For $\Omega\in R_\ell(\A_n^A)$ and $\bm x=(x_1,\ldots,x_n)\in\Omega$, let $\pi$ be the associated permutation and $\lambda(D_{\pi,1}(\bm x))=(B_1,\ldots,B_\ell)$, where $B_k=\{i_{k,1},\ldots,i_{k,n_k}\}$ with $i_{k,1}<\cdots<i_{k,n_k}$ for $1\le k\le \ell$. Define
\[
\bm x^{(k)}=(x_{i_{k,1}},\ldots,x_{i_{k,n_k}}),\; k=1,\ldots,\ell.
\]
For each $k$, let $\Omega_{k}$ be the region of $\A_{n_k}^A$ containing $\bm x^{(k)}$. We then set 
\[
\phi_\ell(\Omega)=((B_1,\ldots,B_\ell),\Omega_1,\ldots,\Omega_\ell). 
\]
\begin{theorem}
	For every integer $\ell$ with $1\le\ell\le n$, the map $\phi_{\ell}$ defined
	in~\eqref{semiorder-map} is bijective.
\end{theorem}
\begin{proof}
For any region $\Omega$ of $\A_n^A$,
Lemma~\ref{partition-independent} shows that the ordered partition in
$\phi_\ell(\Omega)$ is independent of the chosen point. The restriction to a
block $B_k$ has a prime first Dyck path, so the corresponding region
$\Omega_k$ has level one. Thus $\phi_\ell$ is well defined.

Suppose two regions have the same image. Their sign data agree for pairs of
indices in a common block because the corresponding restrictions agree. For
$i\in B_p$ and $j\in B_q$ with $p<q$, both regions satisfy
$x_i-x_j>a_1$; hence their sign data also agree on every cross-block pair.
The two regions are therefore equal, proving injectivity.

Conversely, fix an ordered set partition $(B_1,\ldots,B_\ell)$ and level-one
regions $\Omega_k$ of the induced arrangements on the blocks. After
labeling its coordinates by $B_k$, choose
${\bm x}^{(k)}=(x_i^{(k)})_{i\in B_k}\in\Omega_k$. Set $T_\ell=0$ and choose
$T_{\ell-1},\ldots,T_1$ recursively so that
\[
\min_{i\in B_p}(x_i^{(p)}+T_p)
>
\max_{\substack{q>p\\j\in B_q}}(x_j^{(q)}+T_q)+a_1
\qquad(1\le p<\ell).
\]
Translate all coordinates in $B_k$ by $T_k$.
The resulting point lies in a region $\Omega$ whose ordered prime-component
partition is $(B_1,\ldots,B_\ell)$ and whose restriction to $B_k$ is
$\Omega_k$. Hence $\phi_\ell(\Omega)=((B_1,\ldots,B_\ell),\Omega_1,\ldots,
\Omega_\ell)$, proving surjectivity.
\end{proof}

\section{Characteristic polynomials in the binomial basis}\label{sec-3}

We now identify the coefficients of the characteristic polynomial in the
binomial basis. With the intersection-poset convention fixed in the
introduction,
\[
\chi(\A,t)=\sum_{X\in L(\A)}\mu(\mathbb{R}^n,X)\,t^{\dim X}.
\]

Stanley~\cite{Stanley1996} introduced the notion of an \emph{exponential
sequence of arrangements} (ESA): a sequence
$\mathfrak{A}=(\mathcal{A}_1,\mathcal{A}_2,\ldots)$ of real arrangements
such that
\begin{enumerate}
    \item[(a)] $\mathcal{A}_n$ is an arrangement in $\mathbb{R}^n$;
    \item[(b)] every hyperplane $H\in\mathcal{A}_n$ is parallel to some $x_i-x_j=0$ of the braid arrangement $\mathcal{B}_n$;
    \item[(c)] for any $k$-subset $S\subseteq[n]$, the restricted arrangement
    \[
    \mathcal{A}_n^S=\{H\in\mathcal{A}_n: H \text{ is parallel to } x_i-x_j=0\text{ for some }i,j\in S\}
    \]
    satisfies $L(\mathcal{A}_n^S)\cong L(\mathcal{A}_k)$.
\end{enumerate}
Braid arrangements, Shi arrangements, and symmetric deformations
$(\A_1^A,\A_2^A,\ldots)$ for finite $A\subseteq\mathbb R$ are ESAs; see
\cite[Section~5.3]{Stanley2007}. We use the following result of Stanley.
The face-level analogue for ESAs, together with a Whitney-polynomial
extension, is developed in~\cite{ChenFuLiangWang2026}.

\begin{theorem}[{\cite[Theorem~1.2]{Stanley1996}}]\label{ESA-thm}
Let $\mathfrak{A}=(\A_1,\A_2,\ldots)$ be an ESA. Then
\[
1+\sum_{n\ge 1}\chi(\A_n,t)\frac{x^n}{n!}=\Bigl(1+\sum_{n\ge 1}(-1)^n r(\A_n)\frac{x^n}{n!}\Bigr)^{\!-t}.
\]
\end{theorem}

\begin{proof}[Proof of Theorem~\ref{main-2}]
	Fix $A$ and apply Theorem~\ref{ESA-thm} to
	$\mathfrak{A}^A=(\A_1^A,\A_2^A,\ldots)$; to lighten notation, write
	$\A_n=\A_n^A$. Then
	\begin{equation}\label{eq-1}
	1+\sum_{n\ge1}\chi(\A_n,t)\frac{x^n}{n!}=\Big(1+\sum_{n\ge 1}(-1)^{n}r(\A_n)\frac{x^n}{n!}\Big)^{-t}.
	\end{equation}
	Setting $t=1$ gives
	\begin{equation}\label{eq-2}
		1+\sum_{n \geq 1}\chi({\A_n},1)\frac{x^{n}}{n!}=\left(1+\sum\limits_{n \geq1}(-1)^{n}r(\mathcal{A}_{n})\frac{x^{n}}{n!}\right)^{-1}.
	\end{equation}
	Every $\A_n$ has the diagonal line
	$L=\mathbb R(1,\ldots,1)$ as its translation space. Its
	essentialization in $\mathbb R^n/L$ has rank $n-1$, and its bounded
	regions are precisely the level-$1$ regions of $\A_n$. Zaslavsky's
	formula~\cite[Theorem~C]{Zaslavsky1975} therefore gives
	$\chi(\A_n,1)=(-1)^{n-1}r_1(\A_n)$. Substituting into
	\eqref{eq-2} yields
	\begin{align}\label{eq-3}
		1+\sum_{n \geq 1} (-1)^{n-1}r_1(\mathcal{A}_{n})\frac{x^{n}}{n!}=\left(1+\sum\limits_{n \geq1}(-1)^{n}r(\mathcal{A}_{n})\frac{x^{n}}{n!}\right)^{-1}.
	\end{align}
	Combining \eqref{eq-1}, \eqref{eq-3}, and \eqref{def-function} gives
	\begin{align*}
		1+\sum_{n \geq 1}\chi({\A_n},t)\frac{x^{n}}{n!}
		&=\left( 1+\sum_{n \geq 1} (-1)^{n-1}r_1(\mathcal{A}_{n})\frac{x^{n}}{n!}\right)^t\\[6pt]
		&=\left(1-F_1(\mathfrak{A}^A,-x)\right)^t\\[6pt]
		&=1+\sum_{\ell \geq 1}(-1)^{\ell}\binom{t}{\ell}F_{1}(\mathfrak{A}^A,-x)^\ell.
	\end{align*}
	By Theorem~\ref{main-1}, $F_{1}(\mathfrak{A}^A,-x)^\ell = F_{\ell}(\mathfrak{A}^A,-x)$, so
	\[1+\sum_{n \geq 1}\chi({\A_n},t)\frac{x^{n}}{n!}= 1+\sum_{\ell \geq 1}(-1)^{\ell}\binom{t}{\ell}F_{\ell}(\mathfrak{A}^A,-x).\]
	Substituting \eqref{def-function} and comparing coefficients of $\frac{x^n}{n!}$ on both sides completes the proof.
\end{proof}

\section{Stirling convolution for Catalan- and semiorder-type arrangements}\label{sec-4}
In this section we take $A=\{a_1,\ldots,a_m,0\}$ with $a_1>\cdots>a_m>0$ and $A^*=A\setminus\{0\}$, and prove Theorem~\ref{main-3} via the Stirling convolution~\eqref{Stirling-convolution}. Note that
\begin{align*}
    F_\ell\!\left(\mathfrak{A}^{A^*},\log\frac{1}{1-x}\right)
    &=\sum_{k\ge 1}r_\ell(\A_k^{A^*})\frac{1}{k!}
    \left(\log\frac{1}{1-x}\right)^k
     =\sum_{k\ge 1}r_\ell(\A_k^{A^*})\sum_{n\ge 1}c(n,k)\frac{x^n}{n!}\\
    &=\sum_{n\ge 1}\frac{x^n}{n!}\sum_{k=1}^n c(n,k)\,r_\ell(\A_k^{A^*}).
\end{align*}
Since $F_\ell(\mathfrak{A}^A,1-e^{-x})=F_\ell(\mathfrak{A}^{A^*},x)$
is equivalent to
\[
F_\ell(\mathfrak{A}^A,x)
=F_\ell\!\left(\mathfrak{A}^{A^*},\log\frac{1}{1-x}\right),
\]
it suffices to establish~\eqref{Stirling-convolution}. We do so with the
labeled-Dyck-path model.
For ${\bm x}=(x_1,\ldots,x_n)$ in the complement of
$\A_n^{A^*}$, let $\pi$ be its associated permutation, as defined in
Section~\ref{sec-2}. Associate to
${\bm x}$ the $m$ sign matrices
\begin{equation*}
M_k({\bm x})=\bigl(\operatorname{sgn}(x_{\pi(i)}-x_{\pi(j)}-a_k)\bigr)_{i,j\in[n]},\quad 1\le k\le m,
\end{equation*}
which for $k=1$ agrees with the definitions~\eqref{catalan-M1} and~\eqref{semiorder-M1}. Since $x_{\pi(1)}\ge\cdots\ge x_{\pi(n)}$, all $+$ entries lie above the diagonal and cluster in the northeast of each $M_k({\bm x})$. The boundary between $+$ and $-$ in $M_k({\bm x})$ is therefore a Dyck path $D_{\pi,k}({\bm x})$ labeled by $\pi$, yielding an $m$-tuple
\begin{align}\label{tuple-Dyck}
D_\pi({\bm x})=\bigl(D_{\pi,1}({\bm x}),\ldots,D_{\pi,m}({\bm x})\bigr).
\end{align}

\begin{exam}\label{ex3}
	Use the region $\Omega$ and point
	${\bm x}=(3.0,2.2,4.5,8.0,8.0)$ from Example~\ref{ex2}. The resulting
	$3$-tuple of labeled Dyck paths is shown in Figure~\ref{fig-tuple}.
	\begin{figure}[htbp]
		\centering  
		\tikzset{global scale/.style={
				scale=#1,
				every node/.append style={scale=#1}
			}
		}  
		\begin{tikzpicture}[global scale = 1]
			
			\begin{scope}[scale=0.75, shift={(0,6.5)}]
				\draw (0,0) grid (5,5);
				\draw[red,line width=0.5mm] (0,5) -- (2,5) -- (2,4) -- (2,3) -- (3,3) -- (3,2) -- (5,2) -- (5,0);
				\node at (-0.2,-0.2) {$O$};
				\draw[->, thick] (0,0) -- (5.5,0) node[below] {$x$};
				\draw[->, thick] (0,0) -- (0,5.5) node[above left] {$y$};
				\node at (0.5,5.3) {4};
				\node at (1.5,5.3) {5}; 
				\node at (1.75,4.5) {4};
				\node at (2.5,2.65) {3};
				\node at (1.75,3.5) {5};
				\node at (3.2,2.5) {3};
				\node at (3.5,1.65) {1}; 
				\node at (5.25,1.5) {1};
				\node at (4.5,1.65) {2};
				\node at (5.25,0.5) {2}; 
				\node at (2.5,-0.6) {$D_{\pi,3}(\bm x)$};
			\end{scope}
			
			\begin{scope}[scale=0.75, shift={(-6.5,6.5)}]
				\draw (0,0) grid (5,5);
				\draw[red,line width=0.5mm] (0,5) -- (2,5) -- (2,3) -- (4,3) -- (4,2) -- (5,2) -- (5,0);
				\node at (-0.2,-0.2) {$O$};
				\draw[->, thick] (0,0) -- (5.5,0) node[below] {$x$};
				\draw[->, thick] (0,0) -- (0,5.5) node[above left] {$y$};
				\node at (0.5,5.3) {4};
				\node at (1.5,5.3) {5}; 
				\node at (1.75,4.5) {4};
				\node at (2.5,2.65) {3};
				\node at (1.75,3.5) {5};
				\node at (4.2,2.5) {3};
				\node at (3.5,2.65) {1}; 
				\node at (5.2,1.5) {1};
				\node at (4.5,1.65) {2};
				\node at (5.3,0.5) {2}; 
				\node at (2.5,-0.6) {$D_{\pi,2}(\bm x)$};
			\end{scope}
			
			\begin{scope}[scale=0.75, shift={(-13,6.5)}]
				\draw (0,0) grid (5,5);
				\draw[red,line width=0.5mm] (0,5) -- (2,5) -- (2,3) -- (5,3) -- (5,0);
				\node at (-0.2,-0.2) {$O$};
				\draw[->, thick] (0,0) -- (5.5,0) node[below] {$x$};
				\draw[->, thick] (0,0) -- (0,5.5) node[above left] {$y$};
				\node at (0.5,5.3) {4};
				\node at (1.5,5.3) {5}; 
				\node at (1.75,4.5) {4};
				\node at (2.5,2.65) {3};
				\node at (1.75,3.5) {5};
				\node at (5.3,2.5) {3};
				\node at (3.5,2.65) {1}; 
				\node at (5.3,1.5) {1};
				\node at (4.5,2.65) {2};
				\node at (5.3,0.5) {2}; 
				\node at (2.5,-0.6) {$D_{\pi,1}(\bm x)$};
			\end{scope}  
		\end{tikzpicture}
		\caption{$3$-tuple $D_{45312}(\bm x)$ of labeled Dyck paths}
		\label{fig-tuple}
	\end{figure}
\end{exam}
If $\Delta$ is a region of $\A_n^A$, then
$\mathcal{B}_n\subseteq\A_n^A$ ensures that the associated permutation is
the same for all points in $\Delta$. Moreover,
$\operatorname{sgn}(x_i-x_j-a_k)$ is independent of ${\bm x}\in\Delta$
for all $i,j\in[n]$ and $1\le k\le m$. Hence we can write
\begin{equation}\label{tuple-Dyck1}
D_\pi(\Delta)=\bigl(D_{\pi,1}(\Delta),\ldots,D_{\pi,m}(\Delta)\bigr),
\end{equation}
where
$\Delta\subseteq C_\pi=\{{\bm x}\in\mathbb R^n:
x_{\pi(1)}>\cdots>x_{\pi(n)}\}$.

Now let $\Omega$ be a region of $\A_n^{A^*}$. As in
Lemma~\ref{matrix-independent}, every $M_k({\bm x})$ has the same
unlabeled Dyck path throughout $\Omega$, but its labels may change. For
example, the point ${\bm y}=(3.0,2.2,4.5,8.0,8.1)$ belongs to the region
in Example~\ref{ex3} and yields $D_{54312}({\bm y})$: its three unlabeled
paths agree with those of $D_{45312}({\bm x})$, but the permutations
differ. This ambiguity occurs because $\A_n^{A^*}$ does not contain the
braid arrangement.

We use the following word-partition notation to describe the valid
labelings. For a permutation word $\pi(1)\cdots\pi(n)$, a
\emph{contiguous partition} of $\pi$ is a tuple of contiguous subwords
\[
p=\bigl(\pi(1)\cdots\pi(i_1),\;\pi(i_1+1)\cdots\pi(i_2),\;\ldots,\;\pi(i_{s-1}+1)\cdots\pi(i_s)\bigr)
\]
with $1\le i_1<\cdots<i_s=n$. Each subword is a \emph{block} of $p$, and
$p$ has $s$ blocks. We write $p$ compactly as
\[
p = \pi(1)\cdots\pi(i_1)\mid\pi(i_1+1)\cdots\pi(i_2)\mid\cdots\mid\pi(i_{s-1}+1)\cdots\pi(i_s).
\]
The set of all contiguous partitions of $\pi$ is denoted
$\mathcal{P}(\pi)$. Under refinement it is a poset: $p_1\le p_2$ when
each block of $p_1$ is contained in a block of $p_2$. The meet
$p_1\wedge p_2$ is their common refinement, obtained by retaining every
cut appearing in either partition. For example, if $\pi=76548213$,
$p_1=76|548|213$ and $p_2=765|4821|3$, then
$p_1\wedge p_2=76|5|48|21|3$. Partitions $p$ of $\pi$ and $q$ of
$\sigma$ are \emph{equivalent} if corresponding blocks contain the same
labels.

A labeled Dyck path $D_\pi$ determines a contiguous partition
$p(D_\pi)$: labels belong to the same block when their East steps have the
same height and their South steps have the same horizontal position. This
is an equivalence relation. Moreover, both step coordinates vary
monotonically along the permutation word, so every equivalence class is
contiguous. For instance, $D_{3471526}$ in Figure~\ref{fig-partition}
yields $p(D_{3471526})=347|15|26$.

\begin{figure}[htbp]
	\centering  
	\tikzset{global scale/.style={
			scale=#1,
			every node/.append style={scale=#1}
		}
	}  
	\begin{tikzpicture}[global scale = 1]
		\begin{scope}[scale=0.7, shift={(0,6.5)}]
			\draw (0,0) grid (7,7);
			\draw[red, line width=0.5mm] (0,7) -- (5,7) -- (5,4) -- (7,4) -- (7,0);
			\draw[->, thick] (0,0) -- (7.7,0) node[below] {$x$};
			\draw[->, thick] (0,0) -- (0,7.7) node[left] {$y$};
			\node at (0.5,7.35) {3};
			\node at (1.5,7.35) {4};
			\node at (2.5,7.35) {7};
			\node at (3.5,7.35) {1}; 
			\node at (4.7,6.5) {3};
			\node at (4.7,5.5) {4};
			\node at (4.5,7.35) {5};
			\node at (5.5,3.65) {2};
			\node at (4.65,4.5) {7};
			\node at (7.3,3.5) {1};
			\node at (7.3,1.5) {2};
			\node at (7.3,2.5) {5};
			\node at (6.5,3.65) {6}; 
			\node at (7.3,0.5) {6};
		\end{scope}
	\end{tikzpicture}
	\caption{The labeled Dyck path $D_{3471526}$ and its induced partition
	$347\mid15\mid26$}
	\label{fig-partition}
\end{figure}
\noindent For a region \(\Omega\) of $\A_n^{A^*}$ and any \({\bm x} \in \Omega\) with associated permutation $\pi$, we obtain an \(m\)-tuple of labeled Dyck paths as in \eqref{tuple-Dyck}
\begin{align*}
	D_\pi(\bm x) = (D_{\pi,1}(\bm x), D_{\pi,2}(\bm x),\ldots, D_{\pi,m}(\bm x)).
\end{align*}
Using the above construction, we extract $m$ partitions $p(D_{\pi,1}(\bm x)),\ldots,p(D_{\pi,m}(\bm x))$ of $\pi$, and define their meet in $\mathcal{P}(\pi)$ by
\begin{align}\label{partition-auto}
	p_{\pi}(\bm x):= \bigwedge_{k=1}^m p(D_{\pi,k}(\bm x)).
\end{align}
For a region $\Omega$ of $\A_n^{A^*}$ and $i\in [n]$, pick $\bm x=(x_1,\ldots,x_n)\in\Omega$ and define
\begin{align*}
	\Lambda_k(\Omega,i)=\{j\in [n]\colon x_i-x_j-a_k>0 \}, \quad 
	V_k(\Omega,i)=\{j\in [n]\colon x_j-x_i-a_k>0\}.
\end{align*}
The signs of the defining affine forms are constant on $\Omega$, so these
sets are independent of the choice of $\bm x$. We call $i$ and $j$
\emph{threshold-equivalent} in $\Omega$ if
$\Lambda_k(\Omega,i)=\Lambda_k(\Omega,j)$ and
$V_k(\Omega,i)=V_k(\Omega,j)$ for all $1\le k\le m$. The following lemma
characterizes precisely which permutations can label the Dyck path tuple of
a semiorder-type region.

\begin{lemma}\label{unify}
Let $\Omega$ be a region of $\A_n^{A^*}$ and let ${\bm x}\in\Omega$ have associated permutation $\pi$.
\begin{enumerate}
    \item[(1)] For any ${\bm y}\in\Omega$ with associated permutation $\sigma$, the partitions $p_\pi({\bm x})$ and $p_\sigma({\bm y})$ are equivalent.
    \item[(2)] Conversely, if $\sigma\in\mathfrak{S}_n$ and there exists a partition $q_\sigma$ of $\sigma$ equivalent to $p_\pi({\bm x})$, then some ${\bm y}\in\Omega$ has associated permutation $\sigma$ and satisfies $p_\sigma({\bm y})=q_\sigma$.
\end{enumerate}
\end{lemma}

\begin{proof}
The proof rests on the following claim.
\begin{claim}
For $i,j\in[n]$, the labels $i$ and $j$ lie in the same block of
$p_\pi({\bm x})$ if and only if they are threshold-equivalent in $\Omega$.
\end{claim}
In the $k$-th labeled Dyck path $D_{\pi,k}({\bm x})$, the East step
labeled $i$ has $y$-coordinate
$n-|V_k(\Omega,i)|$, and the South step labeled $i$ has $x$-coordinate
$n-|\Lambda_k(\Omega,i)|$. Hence the East steps labeled $i,j$ have the
same height and the South steps labeled $i,j$ have the same horizontal
position exactly when
\[
|\Lambda_k(\Omega,i)|=|\Lambda_k(\Omega,j)|,
\qquad
|V_k(\Omega,i)|=|V_k(\Omega,j)|.
\]
At any fixed point of $\Omega$, a set $\Lambda_k(\Omega,i)$ of cardinality
$s$ consists of the labels of the $s$ smallest coordinates: each of its
coordinates is strictly smaller than every coordinate outside the set.
Thus its cardinality determines the set. Similarly, a set
$V_k(\Omega,i)$ of cardinality $s$ consists of the labels of the $s$
largest coordinates. Equal cardinalities therefore imply equality of the
corresponding sets, proving the claim.

\noindent\textit{(1)} Since $\Lambda_k(\Omega,i)$ and $V_k(\Omega,i)$
depend only on $\Omega$, the claim gives the same threshold-equivalence
classes at every point of $\Omega$. Their order is also fixed. Indeed, if
two distinct classes reversed order between two points, the segment
joining those points in the convex region $\Omega$ would contain a point
where one coordinate from each class is equal. At that point their
$\Lambda_k$- and $V_k$-sets would agree for every $k$, contradicting that
the classes are distinct. Thus $p_\pi({\bm x})$ and
$p_\sigma({\bm y})$ have corresponding blocks with the same labels and
are equivalent.

\noindent\textit{(2)} Write $p_\pi({\bm x})=B_1|\cdots|B_s$ and
$q_\sigma=C_1|\cdots|C_s$, and let $\sigma=\omega\circ\pi$. Equivalence
means that $B_p$ and $C_p$ contain the same labels for every $p$; hence
$\omega$ maps each threshold-equivalence class $B_p$ to itself. Set
${\bm y}=\omega\cdot{\bm x}$, so $y_i=x_{\omega^{-1}(i)}$. Then $\sigma$
is the associated permutation of ${\bm y}$.

It remains to show that ${\bm y}\in\Omega$. For $i,j\in[n]$ and
$1\le k\le m$, threshold equivalence gives
\[
\begin{aligned}
x_i-x_j>a_k
&\Longleftrightarrow j\in\Lambda_k(\Omega,i)\\
&\Longleftrightarrow j\in\Lambda_k(\Omega,\omega(i))\\
&\Longleftrightarrow \omega(i)\in V_k(\Omega,j)\\
&\Longleftrightarrow \omega(i)\in V_k(\Omega,\omega(j))\\
&\Longleftrightarrow x_{\omega(i)}-x_{\omega(j)}>a_k.
\end{aligned}
\]
No defining form vanishes on $\Omega$, so all corresponding signs agree.
This proves ${\bm y}\in\Omega$. The claim, together with the associated
permutation $\sigma$, now gives $p_\sigma({\bm y})=q_\sigma$.
\end{proof}

\begin{exam}
To illustrate Lemma~\ref{unify}, let $\Omega$ be the region of
$\A_7^{A^*}$ containing
${\bm x}=(4.7,1.6,4.75,4.71,6.6,2.5,1.61)$, with associated permutation
$\pi=5341672$. One computes $p_{\pi}({\bm x})=5|341|6|72$
from~\eqref{partition-auto}. The partition
$q_{\sigma}=5|431|6|27$ of $\sigma=5431627$ is equivalent to
$p_{\pi}({\bm x})$, and $\sigma=\omega\circ\pi$ with
$\omega=1743562$. Lemma~\ref{unify} gives
${\bm y}=(4.7,1.61,4.71,4.75,6.6,2.5,1.6)\in\Omega$ and
$p_{\sigma}({\bm y})=q_{\sigma}$.
\end{exam}

For any region $\Omega$ of $\A_n^{A^*}$, part~(1) of Lemma~\ref{unify} shows that the partitions~\eqref{partition-auto} are equivalent for all points in $\Omega$. By part~(2), we can select a point
\begin{equation}\label{point-choose}
{\bm x}=(x_1,\ldots,x_n)\in\Omega
\end{equation}
whose associated permutation $\pi$ satisfies $p_\pi({\bm x})=B_1|\cdots|B_s$ with each block $B_j$ an increasing consecutive subpermutation. We then define $D_\pi(\Omega)$ as the $m$-tuple $D_\pi({\bm x})$ from~\eqref{tuple-Dyck} and set $p_\pi(\Omega)=p_\pi({\bm x})$.

We briefly recall the cycle representation of permutations; see
\cite[Section~1.3]{Stanley2012}. A permutation $\omega$ is a product of
disjoint cycles. A cycle of length $d$, based at $x$, is written
$(x,\omega(x),\ldots,\omega^{d-1}(x))$. In the \emph{standard
representation}, each cycle begins with its largest element, and the
cycles are ordered increasingly by their largest elements. Removing the
parentheses yields a word $\widehat\omega$; the map
$\omega\mapsto\widehat\omega$ is the classical \emph{fundamental
bijection} of $\mathfrak{S}_n$.

\begin{proof}[Proof of Theorem~\ref{main-3}]
We give a bijective proof. Let \(\operatorname{Cyc}_n(k)\) denote the set of permutations of $[n]$ with exactly \(k\) cycles. We construct a bijection 
\[
\Phi\colon \bigsqcup_{k=1}^n \operatorname{Cyc}_n(k)\times R_\ell(\A_k^{A^*}) \longrightarrow R_\ell(\A_n^A)
\]
as follows. Given \((\omega, \Omega)\) with $\omega\in\operatorname{Cyc}_n(k)$ and $\Omega\in R_\ell(\A_k^{A^*})$, pick $\bm x=(x_1,\ldots,x_k)$ as in \eqref{point-choose} with associated permutation $\pi\in\mathfrak{S}_k$. Let $p_\omega = C_1 | C_2 | \cdots | C_k$ be the partition of $\widehat\omega$ induced by the standard cycle representation of $\omega$, and write each cycle $C_j=i_{j,1}i_{j,2}\cdots i_{j,n_j}$. For $1\le j\le k$, set
\begin{align}\label{con-Phi}
y_{i_{j,1}}&=x_j+\Bigl(\sum_{i=j}^{k}n_i-1\Bigr)\epsilon,\notag\\
y_{i_{j,2}}&=x_j+\Bigl(\sum_{i=j}^{k}n_i-2\Bigr)\epsilon,\quad\ldots,\quad
y_{i_{j,n_j}}=x_j+\sum_{i=j+1}^{k}n_i\epsilon .
\end{align}
Set
\[
\delta({\bm x})=
\min\Bigl(\{a_m\}\cup
\bigl\{|x_p-x_q|:p\ne q,\ x_p\ne x_q\bigr\}\cup
\bigl\{|x_p-x_q-a_s|:p\ne q,\ 1\le s\le m\bigr\}\Bigr)>0
\]
and choose $\epsilon>0$ with $n\epsilon<\delta({\bm x})$. Then, for
every $a_s$, the inequality $x_p-x_q>a_s$ holds if and only if
$y_{i_{p,s_1}}-y_{i_{q,s_2}}>a_s$ for all
$1\le s_1\le n_p$ and $1\le s_2\le n_q$. This yields
${\bm y}\in\mathbb R^n$, which lies in a unique region $\Delta$ of
$\A_n^A$. Set $\Phi(\omega,\Omega)=\Delta$.

We first verify that this definition is independent of the chosen point
${\bm x}$ and of $\epsilon$. Lemma~\ref{unify} shows that requiring every
block of $p_\pi({\bm x})$ to be increasing determines the associated
permutation $\pi$: the threshold-equivalence classes have a fixed order,
and labels within each class are placed increasingly. If $x_p\ne x_q$,
the additional bound $n\epsilon<|x_p-x_q|$ preserves their order in the
constructed point. If $x_p=x_q$, then $p$ and $q$ are
threshold-equivalent; their order is fixed by the increasing order inside
their common block, and the cumulative offsets in~\eqref{con-Phi} impose
that same order on all coordinates in their cycles. Thus the signs of all
braid forms depend only on $\omega$ and $\Omega$. For a nonzero offset
$a_r$, two coordinates in the same cycle differ in absolute value by less than
$n\epsilon<a_m\le a_r$, so their signs are fixed. For coordinates in
distinct cycles $C_p,C_q$, their sign equals that of
$x_p-x_q-a_r$, which is constant on $\Omega$ by the choice of
$\delta({\bm x})$. Thus every admissible choice produces the same region
$\Delta$, and $\Phi$ is well defined.

The inverse construction needed below is recorded separately.

\begin{lemma}[Recovery lemma]\label{lem:recovery}
Let $\Delta$ be a region of $\A_n^A$. The following procedure determines a
unique pair $(\omega,\Omega)$ with
$\omega\in\operatorname{Cyc}_n(k)$ and
$\Omega\in R(\A_k^{A^*})$ for some $k$, and this procedure is inverse to
$\Phi$.
\end{lemma}

\begin{proof}
The braid hyperplanes belong to $\A_n^A$, so every point of $\Delta$ has
the same associated permutation $\sigma$. Choose ${\bm y}\in\Delta$ and
let $\widetilde\Omega$ be the unique region of $\A_n^{A^*}$ containing
$\Delta$. The claim in Lemma~\ref{unify} identifies the
threshold-equivalence classes of $\widetilde\Omega$ with the blocks of
the contiguous partition
\[
p_\sigma(\Delta)=
\bigwedge_{r=1}^m p(D_{\sigma,r}(\Delta))
=B_1|\cdots|B_s.
\]
This partition depends only on the signs defining $\Delta$. In each word
$B_p$, place a cut immediately before every left-to-right maximum. If
$B_p=C_{p,1}|\cdots|C_{p,i_p}$, each word $C_{p,q}$ begins with its
largest entry, and these initial entries increase with $q$. Regard the
$C_{p,q}$ as cycles and put all resulting cycles into standard order. This
defines a unique permutation $\omega$; write its standard cycle word as
$\widehat\omega=C_1\cdots C_k$ and let $j_t$ be the first entry of $C_t$.
Set
\[
{\bm x}=(y_{j_1},\ldots,y_{j_k})
\]
Let $\pi\in\mathfrak{S}_k$ be the associated permutation of ${\bm x}$, and
let $\Omega$ be the region of $\A_k^{A^*}$ containing ${\bm x}$.

We first check that the construction is independent of ${\bm y}$. The
permutation $\sigma$, the tuple of sign matrices, and hence
$p_\sigma(\Delta)$ and all left-to-right-maximum cuts are constant on
$\Delta$. If $C_u$ and $C_v$ are two recovered cycles, then the sign of
$x_u-x_v-a_r=y_{j_u}-y_{j_v}-a_r$ is also constant on $\Delta$ for every
$r$. Thus the resulting region $\Omega$ is independent of the chosen
point.

We also verify that the threshold-equivalence blocks of $\Omega$ are
exactly the collections of cycle representatives coming from
$B_1,\ldots,B_s$. Representatives of cycles cut from the same $B_p$
remain threshold-equivalent because all labels in $B_p$ have identical
$\Lambda_r$- and $V_r$-sets in $\widetilde\Omega$. Conversely, if $B_p\ne B_q$,
some label $h$ witnesses a difference between a corresponding pair of
$\Lambda_r$- or $V_r$-sets. Let $j_w$ represent the recovered cycle
containing $h$. Since $h$ and $j_w$ lie in the same block of
$p_\sigma(\Delta)$, they are threshold-equivalent in $\widetilde\Omega$;
replacing $h$ by
$j_w$ preserves the witnessing membership relation. The difference is
therefore visible among the representatives, so distinct $B_p$ cannot
merge in $\Omega$.

Within a fixed $B_p$, the initial entries created by the
left-to-right-maximum cuts increase. They are the largest entries of their
cycles, so standard cycle order gives the corresponding cycle indices in
increasing order. Consequently, every block of $p_\pi({\bm x})$ is an
increasing consecutive subpermutation. Thus ${\bm x}$ satisfies the
canonical condition in~\eqref{point-choose}, and
$p_\pi(\Omega)=p_\pi({\bm x})$.

Now start with $(\omega,\Omega)$ and construct $\Delta$ by
\eqref{con-Phi}. Coordinates indexed by a common cycle differ by less than
$n\epsilon<a_m$ and occur in the order prescribed by that cycle. Cycles
whose indices lie in a common block of $p_\pi(\Omega)$ concatenate in the
order prescribed by $\pi$; the infinitesimal offsets in
\eqref{con-Phi} make their initial entries precisely the left-to-right
maxima of the resulting word. Hence the above cuts recover exactly the
cycles of $\omega$. For two distinct cycles $C_u,C_v$, the choice of
$\epsilon$ gives, for every $r$,
\[
\operatorname{sgn}(y_i-y_j-a_r)
=\operatorname{sgn}(x_u-x_v-a_r)
\quad(i\in C_u,\ j\in C_v).
\]
The representatives selected by the recovery procedure differ from the
original $x_u$ by less than $n\epsilon$; the definition of
$\delta({\bm x})$ therefore puts the recovered representative point in the
same region $\Omega$. Thus recovery after $\Phi$ returns
$(\omega,\Omega)$.

Conversely, start from $\Delta$ and apply the recovery procedure. Within
each recovered cycle all coordinate differences have absolute value less
than $a_m$. Indeed, if $y_i-y_j>a_m$, then
$j\in\Lambda_m(\widetilde\Omega,i)$ but
$j\notin\Lambda_m(\widetilde\Omega,j)$, contradicting threshold
equivalence; the case $y_j-y_i>a_m$ is symmetric. Thus replacing the
coordinates by the ordered infinitesimal offsets in
\eqref{con-Phi} preserves their braid order and the sign of every form with
nonzero offset. Between distinct cycles, the preceding displayed identity
preserves every nonzero-offset sign. For the braid forms, unequal recovered
representatives retain their order because $n\epsilon$ is smaller than
their distance. Equal representatives belong to the same
threshold-equivalence class; within that class, the left-to-right-maximum
cuts and standard cycle order agree with the cumulative order of the
offsets in~\eqref{con-Phi}. Thus their braid order is also recovered.
Therefore the reconstructed point has exactly the same signs for all
defining equations of $\A_n^A$ as ${\bm y}$ and lies in $\Delta$. Hence
$\Phi(\omega,\Omega)=\Delta$, completing the proof.
\end{proof}

The recovery lemma makes $\Phi$ bijective. It remains to prove that it
preserves level. Let $\Omega\in R(\A_k^{A^*})$ have level $\ell$, and let
$\omega$ be a permutation of $[n]$ with $k$ cycles. Write
$p_\omega=C_1|C_2|\cdots|C_k$, choose ${\bm x}\in\Omega$ as
in~\eqref{point-choose}, with associated permutation
$\pi\in\mathfrak{S}_k$, define ${\bm y}$ by~\eqref{con-Phi}, and put
$\Delta=\Phi(\omega,\Omega)$. By the
level description for semiorder-type regions, there are $\ell$ subsets
\[
\begin{aligned}
B_1&=\{\pi(1),\ldots,\pi(k_1)\},\\
B_2&=\{\pi(k_1+1),\ldots,\pi(k_2)\},\quad\ldots,\quad
B_\ell=\{\pi(k_{\ell-1}+1),\ldots,\pi(k_\ell)\}.
\end{aligned}
\]
with $0=k_0<k_1<\cdots<k_\ell=k$. Consecutive coordinates in a common block
differ by less than $a_1$, whereas
\[
x_i-x_j>a_1
\qquad(i\in B_p,\ j\in B_q,\ p<q).
\]
Let $C_{B_p}=\bigcup_{j\in B_p}C_j$ for $1\le p\le \ell$. From the construction of $\bm y$, for $i\in C_{B_p}$ and $j\in C_{B_q}$,
\[
\begin{cases}
	|y_i-y_j|<na_1,&\text{if $p=q$};\\
	y_i-y_j>a_1,&\text{if $p<q$}.
\end{cases}
\]
Define $W=\{(z_1,\ldots,z_n)\in\mathbb{R}^n\colon z_i=z_j
\text{ if }i,j\in C_{B_p}\}$, and let $C\subseteq W$ be the
full-dimensional cone whose component values decrease with $p$. The
displayed bounds give $\ell(\Delta)\le\dim W=\ell$. A translate of $C$ is
contained in $\Delta$; applying the scaling argument from
Proposition~\ref{components-1} to any subspace whose bounded neighborhood
contains $\Delta$ gives $\ell\le\ell(\Delta)$. Thus
$\ell(\Delta)=\ell$.
\end{proof}

The following example illustrates the bijection.

\begin{exam}
Let $A=\{0,1,2\}$ and $n=7$. Let $\Omega$ be the region of
$\A_5^{[2]}$ containing
${\bm x}=(4.75,4.7,6.6,2.5,1.6)$, and take
$\omega=(3)(41)(5)(6)(72)\in\mathfrak{S}_7$, which has five cycles. With
$\epsilon=0.01$, formula~\eqref{con-Phi} gives
\[
{\bm y}=(4.74,1.60,4.81,4.75,6.63,2.52,1.61).
\]
This point determines a region $\Delta$ of $\A_7^{[0,2]}$, and
$\Phi(\omega,\Omega)=\Delta$. Figure~\ref{fig:tikz7} illustrates the
construction via labeled Dyck paths.

\begin{figure}[htbp]
\centering 
\tikzset{global scale/.style={
scale=#1,
every node/.append style={scale=#1}
}
} 
\begin{tikzpicture}[global scale = 0.65]

\begin{scope}[scale=0.8, shift={(6.5,5)}]
\draw (0,0) grid (5,5);
\draw[red, line width=0.5mm] (0,5) -- (3,5) -- (3,2) -- (5,2) -- (5,0);
\draw[->, thick] (0,0) -- (5.5,0) node[below] {$x$};
\draw[->, thick] (0,0) -- (0,5.5) node[left] {$y$};
\node at (0.5,5.3) {(5)};
\node at (1.5,5.3) {(3)};
\node at (2.5,5.3) {(41)}; 
\node at (3.5,1.65) {(6)};
\node at (4.5,1.65) {(72)};
\node at (2.6,4.5) {(5)};
\node at (2.6,3.5) {(3)};
\node at (2.5,2.5) {(41)}; 
\node at (5.35,1.5) {(6)};
\node at (5.5,0.5) {(72)};
\node at (6,-1.5) {$D_{31245}(\Omega)$ labeled by $\omega=(3)(41)(5)(6)(72)$};
\end{scope}
\begin{scope}[scale=0.8, shift={(13,5)}]
\draw (0,0) grid (5,5);
\draw[red,line width=0.5mm] (0,5) -- (1,5) -- (1,4) -- (3,4) -- (3,2) -- (5,2) -- (5,0);
\draw[->, thick] (0,0) -- (5.5,0) node[below] {$x$};
\draw[->, thick] (0,0) -- (0,5.5) node[left] {$y$};
\node at (0.5,5.3) {(5)};
\node at (1.5,4.5) {(3)};
\node at (2.5,4.5) {(41)}; 
\node at (3.5,1.65) {(6)};
\node at (4.5,1.65) {(72)};
\node at (0.6,4.5) {(5)};
\node at (2.6,3.4) {(3)};
\node at (2.5,2.5) {(41)}; 
\node at (5.5,1.5) {(6)};
\node at (5.5,0.5) {(72)};
\end{scope}

\begin{scope}[scale=0.6, shift={(28,6.5)}]
	\draw[->, thick] (-3,3.5) -- (-1.2,3.5);
	\node at (-2.1,4) {$\Phi$};
	\draw (0,0) grid (7,7);
	\draw[->, thick] (0,0) -- (7.7,0) node[below] {$x$};
	\draw[->, thick] (0,0) -- (0,7.7) node[left] {$y$};
	\draw[red, line width=0.5mm] (0,7) -- (4,7) -- (4,3) -- (7,3) -- (7,0);
	\node at (0.5,7.4) {5};
	\node at (1.5,7.4) {3};
	\node at (2.5,7.4) {4};
	\node at (3.5,7.4) {1}; 
	\node at (3.7,6.5) {5};
	\node at (3.7,5.5) {3};
	\node at (3.7,4.5) {4};
	\node at (3.7,3.5) {1};
	\node at (4.5,2.55) {6};
	\node at (5.5,2.55) {7};
	\node at (7.3,1.5) {7};
	\node at (7.3,2.5) {6};
	\node at (6.5,2.55) {2}; 
	\node at (7.3,0.5) {2};
	\node at (8,-1.8) {$D_{5341672}(\Delta)$};
\end{scope}
\begin{scope}[scale=0.6, shift={(37,6.5)}]
	\draw (0,0) grid (7,7);
	\draw[->, thick] (0,0) -- (7.7,0) node[below] {$x$};
	\draw[->, thick] (0,0) -- (0,7.7) node[left] {$y$};
	\draw[red, line width=0.5mm] (0,7) --(1,7) -- (1,6) -- (4,6) -- (4,3) -- (6,3) -- (7,3) -- (7,0);
	\node at (0.5,7.3) {5};
	\node at (1.5,6.4) {3};
	\node at (2.5,6.4) {4};
	\node at (3.5,6.4) {1}; 
	\node at (0.7,6.5) {5};
	\node at (3.7,5.35) {3};
	\node at (3.7,4.5) {4};
	\node at (3.7,3.5) {1};
	\node at (4.5,2.55) {6};
	\node at (5.5,2.55) {7};
	\node at (7.3,1.5) {7};
	\node at (7.3,2.5) {6};
	\node at (6.5,2.55) {2}; 
	\node at (7.3,0.5) {2};
\end{scope}
\end{tikzpicture}
\caption{An illustration of the bijection $\Phi$ via labeled Dyck paths}
\label{fig:tikz7}
\end{figure}
\end{exam}

\section{Level enumeration for m-Catalan arrangements}\label{sec-6}

When $A=[0,m]$, the resulting arrangement $\A_n^{[0,m]}$ is the
$m$-Catalan arrangement. We prove Theorem~\ref{main-4}. We use the formal
power series
$
B_s(x)=\sum_{q\ge0}\frac{1}{sq+1}\binom{sq+1}{q}x^q,
$
with notation as in~\cite[pp.~201, 364]{GRK-book}. For every positive
integer $s$ and every real $r$, Lagrange inversion gives
\begin{align}\label{eq:Cat-2}
B_s(x)^r
=1+\sum_{q\ge1}\frac{r}{q}\binom{sq+r-1}{q-1}x^q.
\end{align}

\begin{proof}[Proof of Theorem~\ref{main-4}]
The level-$1$ regions of $\A_n^{[0,m]}$ are precisely its relatively
bounded regions. Combining Zaslavsky's formula
\cite[Theorem~C]{Zaslavsky1975} with Athanasiadis's characteristic
polynomial for the $m$-Catalan arrangement
\cite[Theorem~5.1, p.~221]{Athanasiadis1996}, we obtain
\[
r_{e,1}(\A_n^{[0,m]})
=\frac{1}{(m+1)n-1}\binom{(m+1)n-1}{n}.
\]
From \eqref{def-function} and \eqref{eq:Cat-2}, we have
\[
F_1(\mathfrak{A}^{[0,m]}, x) = 1-B_{m+1}(x)^{-1}.
\]
Substituting this into Theorem~\ref{main-1},
\begin{align*}
F_\ell(\mathfrak{A}^{[0,m]},x)
&=\left(1-B_{m+1}(x)^{-1}\right)^\ell\\
&=\sum_{i=0}^\ell(-1)^i\binom{\ell}{i}B_{m+1}(x)^{-i}\\
&=\sum_{n\ge1}\frac{1}{n}
\left(\sum_{i=1}^\ell(-1)^{i-1}i\binom{\ell}{i}
\binom{(m+1)n-i-1}{n-1}\right)x^n,
\end{align*}
where the last equality applies~\eqref{eq:Cat-2} with $r=-i$. Since
$F_1$ has zero constant term, coefficients with $n<\ell$ vanish. Fix
$n\ge\ell$ and compare coefficients to obtain
\begin{align}\label{rel}
r_{e,\ell}(\A_n^{[0,m]})
=\frac{1}{n}\sum_{i=1}^\ell(-1)^{i-1}i\binom{\ell}{i}
\binom{(m+1)n-i-1}{n-1}.
\end{align}
To simplify, we use the binomial identity
$
\sum_{j=0}^a(-1)^j\binom{a}{j}\binom{N-j}{b}
=\binom{N-a}{N-b},
$
which follows by extracting the coefficient of $z^b$ from
\[
\sum_{j=0}^a(-1)^j\binom aj(1+z)^{N-j}
=z^a(1+z)^{N-a}.
\]
Applying this identity to~\eqref{rel},
\begin{align*}
r_{e,\ell}(\A_n^{[0,m]})
&=\frac{\ell}{n}\sum_{j=0}^{\ell-1}(-1)^j
\binom{\ell-1}{j}\binom{(m+1)n-j-2}{n-1}\\
&=\frac{\ell}{n}\binom{(m+1)n-1-\ell}{n-\ell}\\
&=\frac{m\ell}{(m+1)n-\ell}\binom{(m+1)n-\ell}{n-\ell}.
\end{align*}
The last expression is $\operatorname{Ran}_{m+1,m\ell}(n-\ell)$. Multiplying by
$n!$ for the $\mathfrak{S}_n$ action therefore completes the proof.
\end{proof}

\appendix
\section{Tableau realization of the inverse bijection}\label{sec-7}

Recall Stanley's formula relating the number of regions of the $m$-Catalan
arrangement $\A_n^{[0,m]}$ to a Fuss--Catalan number:
\begin{align}\label{eq:region-m-Cat}
\frac{1}{n!}r(\A_n^{[0,m]}) = \operatorname{Ran}_{m+1,1}(n).
\end{align}
A combinatorial proof was given independently by Duarte and Guedes de
Oliveira~\cite{D-G2021} and by Fu, Wang, and Zhu~\cite{F-W2021}. The latter
construct an injection from the regions of $\A_n^{[0,m]}$ in the
fundamental chamber to $m$-Dyck paths of length $(m+1)n$. Krattenthaler
\cite{Krattenthaler1989} showed that the number of such paths is
$\operatorname{Ran}_{m+1,1}(n)$, which is also the number of
fundamental-chamber regions by~\eqref{eq:region-m-Cat}; hence the injection
is a bijection.
Conversely, the labeled-Dyck-path tuple permits an explicit reconstruction
of the region in the fundamental chamber. The reconstruction below is
formulated so that its output has the prescribed column-count data of
Fu--Wang--Zhu~\cite[Theorem~2.6]{F-W2021}.

We now describe the tableau reconstruction. An \emph{$m$-Dyck path}
$D^m$ of length $(m+1)n$ is a lattice path from $(0,mn)$ to $(n,0)$
with $n$ East steps and $mn$ South steps that stays weakly above the line
$y=m(n-x)$. Let $\mathcal D_n^m$ denote the set of these paths. For
$1\le i\le n$, let $h_i'(D^m)$ be the $y$-coordinate at which $D^m$
meets the line $x=i-\frac12$, and define
\[
h_i(D^m)=mn-h_i'(D^m).
\]
The height sequence $h(D^m)=(h_1(D^m),\ldots,h_n(D^m))$ is nondecreasing
and satisfies $0\le h_i(D^m)\le m(i-1)$. Conversely, every nondecreasing
integer sequence with these bounds determines a unique $m$-Dyck path.

\begin{theorem}\label{map:inv}
For integers $m,n\ge1$, the tableau reconstruction algorithm below,
together with a labeling permutation,
defines a bijection from $\mathfrak{S}_n\times\mathcal D_n^m$ to
$R(\A_n^{[0,m]})$. It is inverse to the map of Fu, Wang, and
Zhu~\cite[Theorem~2.6]{F-W2021}. Consequently,
\[
r(\A_n^{[0,m]}) = \frac{n!}{(m+1)n+1} \binom{(m+1)n+1}{n}.
\]
\end{theorem}

Order the positive parameters as $a_k=m-k+1$ for $1\le k\le m$. For a
height sequence $h=(h_1,\ldots,h_n)$, let $\mathcal H(h)$ be the finite set
of integer arrays $H=(h_{k,j})$ satisfying
\begin{equation}\label{eq:height-array}
\begin{gathered}
0=h_{k,1}\le h_{k,2}\le\cdots\le h_{k,n},\qquad
0\le h_{k,j}\le j-1,\\
h_{1,j}\le h_{2,j}\le\cdots\le h_{m,j},\qquad
\sum_{k=1}^m h_{k,j}=h_j.
\end{gathered}
\end{equation}
Associate with $H$ the strict system
\begin{equation}\label{eq:feasible-height-array}
z_1>\cdots>z_n,\qquad
\begin{cases}
z_i-z_j>a_k,&1\le i\le h_{k,j},\\
z_i-z_j<a_k,&h_{k,j}<i<j,
\end{cases}
\quad 1\le k\le m, 1\le j\le n.
\end{equation}
We call $H$ \emph{feasible} if this system has a real solution.

\begin{algorithm}[H]
\SetAlgoLined
\KwIn{An $m$-Dyck path $D^m$ with height sequence $h$}
\KwOut{A tableau $T$ encoding the inverse region in the fundamental chamber}

Form $\mathcal H(h)$ using~\eqref{eq:height-array}\;
Test~\eqref{eq:feasible-height-array} for every $H\in\mathcal H(h)$ and
retain the unique feasible array $H=(h_{k,j})$\;
Initialize an empty $(n-1)\times m$ grid $T$\;

\For{$k\leftarrow1$ \KwTo $m$}{
  \For{$j\leftarrow2$ \KwTo $n$}{
    Insert $h_{k,j}-h_{k,j-1}$ copies of $j$, from top to bottom, into
    column $m-k+1$ of $T$\;
  }
}
\KwRet $(T,H)$
\caption{Tableau reconstruction of the Fu--Wang--Zhu inverse}\label{Al}
\end{algorithm}

\begin{proof}
Fix $\pi\in\mathfrak{S}_n$ and $D^m\in\mathcal D_n^m$.
For a point $\bm z$ in the fundamental chamber, the positive entries in
the $j$th column of its $k$th sign matrix form an initial segment; let its
length be $h_{k,j}$. These lengths satisfy~\eqref{eq:height-array}, and
their sum over $k$ is the Fu--Wang--Zhu height $h_j$. Conversely, a
feasible array and a solution of~\eqref{eq:feasible-height-array}
determine the complete sign-matrix tuple and hence a region.

Theorem~2.6 of Fu--Wang--Zhu states that the column-count map is a
bijection. Thus, for every $m$-Dyck height sequence $h$, exactly one array
in $\mathcal H(h)$ is feasible, so Algorithm~\ref{Al} is well defined.
Column $m-k+1$ of $T$ contains exactly $h_{k,j}$ entries not exceeding
$j$, and therefore $T$ recovers $H$ without loss. The inequalities in
\eqref{eq:height-array} show that these columns form a Young diagram and
are weakly increasing both down columns and across rows. Hence $T$ encodes
the unique fundamental-chamber region mapped to $D^m$. Undoing the
coordinate reordering specified by $\pi$ gives the unique region mapped to
$(\pi,D^m)$. This proves that the construction is inverse to the
Fu--Wang--Zhu map and is a bijection.
\end{proof}

The following example illustrates the algorithm.

\begin{exam}
Consider $m=n=6$, the height sequence
$h=(0,2,4,5,6,12)$, and $\pi=543261$. The feasibility step returns the
array displayed below. The associated tableau, written by rows, is
\[
T=\begin{array}{llllll}
2&2&3&4&6&6\\
3&5&6\\
6\\
6\\
6
\end{array}.
\]

The feasible height array is
\begin{equation*}
\begin{bmatrix}
h_1\\
h_2\\
h_3\\
h_4\\
h_5\\
h_6\\
\end{bmatrix}
=
\begin{bmatrix}
0 & 0 & 0 & 0 & 0 & 1 \\
0 & 0 & 0 & 0 & 0 & 1 \\
0 & 0 & 0 & 1 & 1 & 1 \\
0 & 0 & 1 & 1 & 1 & 2 \\
0 & 1 & 1 & 1 & 2 & 2 \\
0 & 1 & 2 & 2 & 2 & 5 \\
\end{bmatrix},
\end{equation*}
Here the entries not exceeding $j$ in column $6-k+1$ of $T$ recover
$h_{k,j}$. Interpreting the rows as height vectors gives
$D=(D_1,\ldots,D_6)$. Labeling by $\pi$ yields
\[
D_\pi=(D_{\pi,1},D_{\pi,2},D_{\pi,3},D_{\pi,4},D_{\pi,5},D_{\pi,6}).
\]
The point
${\bm x}=(6.1,12.53,13.45,14.51,16.6,12.49)$ lies in a region $\Delta$ of
$\A_6^{[0,6]}$ satisfying $D_\pi(\Delta)=D_\pi$.

\end{exam}

\section*{Acknowledgements}
The second author is supported by the National Natural Science Foundation of
China (Grant No. 12571350) and the Guangdong Basic and Applied Basic Research Foundation
(Grant No. 2025A1515010457). The fourth author is supported by the National Natural Science
Foundation of China (Grant No. 12101613).

\section*{AI Use Statement}

During the preparation of this manuscript, the authors used OpenAI's
ChatGPT and Codex to assist with language editing, computational checks,
and the organization of the exposition. All mathematical results, proofs,
and conclusions presented as contributions in this manuscript are the
authors' original work. The authors reviewed and revised all AI-assisted
output and take full responsibility for the content of the manuscript.

\end{document}